\documentclass[11pt]{article} 
\usepackage[latin1]{inputenc}
\usepackage[T1]{fontenc}
\usepackage[frenchb]{babel}
\usepackage[babel=true,kerning=true]{microtype}
\usepackage{lmodern}
\usepackage{amsmath}
\usepackage{amsfonts} 
\usepackage{amssymb} 
\usepackage{amsthm}
\usepackage{graphicx}
\usepackage{stmaryrd} 
\usepackage[np]{numprint}
\usepackage{url}
\usepackage{vmargin}
\usepackage{pgf,tikz}
\usetikzlibrary{arrows}
\usepackage{mathrsfs}
\theoremstyle{plain}
\newtheorem{theoreme}{Théorème}
\newtheorem{lemme}[theoreme]{Lemme}
\newtheorem{propriete}[theoreme]{Propriété}
\newtheorem{corollaire}[theoreme]{Corollaire}
\theoremstyle{definition}
\newtheorem{definition}[theoreme]{D\'efinition}

\theoremstyle{remark}
\newtheorem{rem}[theoreme]{Remarque}

\newcommand{\ensemblenombre }[1]{\mathbb{#1}}
\newcommand{\N}{\ensemblenombre{N}}
\newcommand{\Z}{\ensemblenombre{Z}}
\newcommand{\Q}{\ensemblenombre{Q}}
\newcommand{\R}{\ensemblenombre{R}}
\newcommand{\C}{\ensemblenombre{C}}

\newcommand{\entiere}[1]{\left\lfloor #1 \right\rfloor}
\newcommand{\intervalleff }[2]{\left[{#1}\mathpunct{};{#2}\right]}
\newcommand{\intervalleof }[2]{\left]{#1}\mathpunct{};{#2}\right]}
\newcommand{\intervallefo }[2]{\left[{#1}\mathpunct{};{#2}\right[}
\newcommand{\intervalleoo }[2]{\left]{#1}\mathpunct{};{#2}\right[}
\newcommand{\diff}{\mathop{}\mathopen{}\mathrm{d}}
\newcommand{\couple}[2]{\left({#1}\mathpunct{};{#2}\right)}
\newcommand{\abs}[1]{\left\lvert#1\right\rvert}

\newcommand*{\EQ}[2]%
{\ensuremath{%
    #1/\!\raisebox{-.65ex}{\ensuremath{{#2}}}}}
\setmarginsrb{2cm}{2cm}{2cm}{2.5cm}{0cm}{0cm}{0cm}{1cm}
\date{} 
\begin{document}

\title{\textbf{Sur l'irrationalité des racines \\ de certaines familles de polynômes}}

\author{Lionel Ponton \\
\small \texttt{lionel.ponton@gmail.com}}

\maketitle

\begin{abstract}
On s'intéresse à l'irrationalité des racines de sept grandes familles de polynômes: les polynômes de Tchebichef, de Legendre, de Laguerre, d'Hermite, de Bessel, de Bernoulli et d'Euler. On établit d'abord, pour la plupart d'entre eux, l'existence de racines réelles avant d'étudier l'irrationalité de celles-ci. Les techniques utilisées reposent le plus souvent sur des propriétés arithmétiques simples des nombres algébriques, certaines étant issues de preuves plus générales ayant permis d'établir l'irréductibilité de certains de ces polynômes.
\end{abstract} 

\bigskip

\bigskip

\section{Introduction}
L'étude de l'irréductibilité sur $\Q$ de certaines familles usuelles de polynômes a connu un essor considérable depuis le début du XX\up{e} siècle. Dans deux articles datant de 1912 et 1913, Hold (\cite{Hol12} et \cite{Hol13}) a établi l'irréductibilité des polynômes de Legendre de degrés particuliers. Quelques années plus tard, Schur \cite{Sch29} a obtenu le premier résultat de portée très générale en établissant, dans un article en deux parties datant de 1929, l'irréductibilité de vastes classes de polynômes incluant les polynômes de Laguerre et les polynômes d'Hermite (à un facteur $X$ près dans le cas d'un degré impair). En 1951, Grosswald \cite{Gro51} a commencé à s'intéresser aux propriétés des polynômes de Bessel avant que Filaseta et Trifonov \cite{FT02} établissent leur irréductibilité en 2002.

Il est intéressant de constater que toutes ces démonstrations s'appuient, pour les premières, sur le critère d'Eisenstein ou, pour les suivantes, sur une généralisation de ce critère datant de 1906: le théorème de Dumas \cite{Dum06}. Ce théorème qui permet, en utilisant les polygones de Newton, de s'assurer qu'un polynôme de $\Z[X]$ n'a pas de facteur d'un degré donné, nécessite, pour obtenir l'irréductibilité de polynômes dans des cas très généraux, d'avoir recours à des résultats profonds sur la répartition des nombres premiers. Ainsi, dans son article de 1929, Schur utilise une généralisation du postulat de Bertrand et, dans la démonstration de l'irréductibilité des polynômes de Bessel, Filaseta et Trifonov utilisent des estimations fines des fonctions $\pi$ et $\theta$ définies respectivement par 
\[\pi(x)=\sum_{p\leqslant x} 1 \qquad \text{et} \qquad \theta(x)=\sum_{p\leqslant x} \ln(p),\]
les sommes portant sur les nombres premiers.

Notre objectif dans ce qui suit est bien plus modeste. L'irréductibilité d'un polynôme sur $\Q$ implique en particulier qu'il n'a pas de facteur affine et donc pas de racine rationnelle. Nous nous proposons d'étudier directement, sans avoir recours à l'irréductibilité, l'irrationalité des racines de sept grandes familles de polynômes: les polynômes de Tchebichef, de Legendre, de Laguerre, d'Hermite, de Bessel, de Bernoulli et d'Euler. Les quatre premières familles font parties d'une famille plus vaste de polynômes appelés, généralement, les polynômes orthogonaux \og classiques \fg{} dont on étudie certaines propriétés dans le paragraphe \ref{para_poly_orthogonaux}. Si on excepte le cas des polynômes de Tchebichef, toutes les démonstrations s'appuient sur des propriétés de la valuation $p-$adique qui sont établies dans le paragraphe \ref{para_valuation_p_adique}. Pour les polynômes de Legendre, de Laguerre, d'Hermite et de Bessel, l'outil central pour obtenir l'irrationalité des racines est la propriété \ref{prop_cri_irra_poly_Newton} qui peut être vue comme un cas particulier d'une propriété plus générale découlant du théorème de Dumas (voir la remarque \ref{rem_poly_Newton}) et que nous établissons directement dans le paragraphe \ref{para_critère_irrationalité}. L'étude de la nature arithmétique des racines des polynômes de Bernoulli et d'Euler repose quant à elle sur le théorème de von Staudt et Clausen (théorème \ref{theo_vS_C}) qui permet d'obtenir une expression des dénominateurs de nombres de Bernoulli (corollaire \ref{coro_vS_C}). 

Nous avons fait le choix, sauf pour les polynômes de Tchebichef et de Legendre, de définir les polynômes étudiés par leurs expressions explicites car les démonstrations d'irrationalité utilisent en général des propriétés arithmétiques des coefficients. Ces polynômes apparaissent souvent, dans certains problèmes mathématiques ou physiques, comme solutions d'équations différentielles ou à travers des relations de récurrence mais nous n'avons pas abordé ces points, sauf s'ils étaient nécessaires pour l'étude des racines. Nous avons, cependant, donné des références bibliographiques traitant ces aspects. 

Dans la plupart des cas, les preuves que nous donnons reprennent des démonstrations existantes ou s'inspirent, en les simplifiant, de démonstrations d'irréductibilité et nous avons systématiquement donné les références des sources utilisées. Lorsqu'aucune source n'est indiquée, le preuve est, à notre connaissance, originale. Même si elles peuvent parfois être un peu techniques (surtout dans le cas des polynômes de Bessel), les démonstrations restent à un niveau élémentaire et ne dépassent pas le cadre du programme de CPGE scientifique tout en permettant d'en aborder un grand nombre de notions: polynômes, arithmétique, espace préhilbertien, intégrales généralisées, séries entières, groupes, équations différentielles...

Pour finir, signalons que ce texte doit beaucoup aux divers articles de M. Filaseta ainsi qu'aux notes du cours \og Math 788F: The Theory  of Irreducible Polynnomials \fg{} qu'il dispense à l'Université de Caroline du Sud.

\section{Résultats préliminaires}
\subsection{Notations et conventions}
Dans toute la suite, on adopte les notations et conventions suivantes.
\begin{enumerate}
\item[$\bullet$] Pour tout réel $x$, $\left\lfloor x\right\rfloor$ désigne la partie entière de $x$ i.e. l'unique entier $k$ tel que $k \leqslant x < k+1$.
\item[$\bullet$] Pour tout nombre premier $p$ et tout entier naturel non nul $n$, on appelle \emph{valuation $p-$adique de $n$}, notée $v_p(n)$, l'unique entier naturel $m$ tel que $p^m$ divise $n$ et $p^{m+1}$ ne divise pas $n$. On convient, de plus, que $v_p(0)=+\infty$.
\item[$\bullet$] On note $\mathbb{P}$ l'ensemble des nombres premiers.
\item[$\bullet$] Pour tous entiers $m$ et $n$ tels que $m\leqslant n$, on note $\llbracket m, n\rrbracket$ l'ensemble des entiers $k$ tels que $m \leqslant k \leqslant n$.
\item[$\bullet$] Si $r$ est un nombre rationnel, \emph{le} numérateur (resp. \emph{le} dénominateur) de $r$ désigne le numérateur (resp. le dénominateur) de l'écriture sous forme de fraction irréductible de $r$.
\item[$\bullet$] Si $P$ est un polynôme de $\R[X]$, on note également $P$ la fonction polynomiale $x \mapsto P(x)$ définie sur $\R$.
\end{enumerate}

\subsection{Polynômes orthogonaux}
\label{para_poly_orthogonaux}

Dans tout ce paragraphe, $I$ désigne un intervalle de la forme $\intervalleoo{a}{b}$ avec $\couple{a}{b}\in\overline{\R}^2$.

\begin{definition} --- On appelle \emph{fonction poids} sur $I$ toute fonction $w: I \to \R_+$ continue et non identiquement nulle sur $I$ telle que, pour tout $P\in\R[X]$, la fonction $t\mapsto P(t)w(t)$ est intégrable sur $I$.
\end{definition}

\begin{propriete} --- Soit $w$ une fonction poids sur $I$. Alors, l'application $\varphi_w$ définie sur $\R[X]^2$ par
\[\forall \couple{P}{Q}\in\R[X]^2 \quad \varphi_w(P,Q):=\int_I P(t)Q(t)w(t) \diff t\]
est un produit scalaire sur $\R[X]$.
\label{prop_prod_scalaire_w}
\end{propriete}

\emph{Preuve}. --- L'application $\varphi_w$ est clairement bilinéaire par linéarité de l'intégrale et symétrique par symétrie du produit réel. De plus, comme $w$ est à valeur positives, pour tout polynôme $P\in\R[X]$, 
\[\varphi_w(P,P)=\int_I P(t)^2w(t) \diff t \geqslant 0\]
avec égalité si $P=0$. Supposons que $P\in\R[X]$ vérifie $\varphi_w(P,P)=0$. Alors, comme $t\mapsto P(t)^2w(t)$ est positive et continue sur $I$, pour tout $t\in I$, $P(t)^2w(t)=0$. Or, $w$ est continue et non identiquement nulle donc il existe un intervalle $J\subset I$ tel que $w(t)>0$ pour tout $t\in J$. Par suite, $P(t)=0$ pour tout $t\in J$ donc $P$ admet une infinité de racines ce qui assure que $P=0$. On conclut que $\varphi_w$ est bien un produit scalaire sur $\R[X]$. \hfill $\square$

\begin{definition} --- Soit $w$ une fonction poids sur $I$. On dit qu'une suite de polynômes $(P_n)\in\R[X]^{\N}$ est \emph{une famille de polynômes $w-$orthogonaux} si, pour tout $n\in\N$, $\deg P_n=n$ et, pour tout $\couple{m}{n}\in\N^2$ tel que $m\neq n$, $\varphi_w(P_m, P_n)=0$ (où $\varphi_w$ est le produit scalaire défini dans la propriété \ref{prop_prod_scalaire_w}).
\end{definition}

\begin{lemme} --- Soit $w$ une fonction poids sur $I$. On munit $\R[X]$ du produit scalaire $\varphi_w$. Alors, une famille $(P_n)$ de polynômes $w-$orthogonaux est une base orthogonale de $\R[X]$ et, pour tout $\couple{m}{n}\in(\N^*)^2$ tel que $n \geqslant m$, $P_n\in\R_{m-1}[X]^{\perp}$.
\label{lem_base_ortho}
\end{lemme}

\emph{Preuve}. --- Soit $(P_n)$ une famille de polynômes $w-$orthogonaux. Comme la famille $(P_n)$ est échelonnée en degré, elle forme une base de $\R[X]$ et elle est, par définition, orthogonale. Soit $\couple{m}{n}\in(\N^*)^2$ tel que $n \geqslant m$. La famille $(P_k)_{k\in\llbracket 0, m-1\rrbracket}$ est une base de $\R_{m-1}[X]$ pour la même raison. Soit $Q\in\R_{m-1}[X]$. Alors, il existe des réels $\lambda_0$, ..., $\lambda_{m-1}$ tels que 
\[Q=\sum_{k=0}^{m-1} \lambda_k P_k.\]
Dès lors, comme $n>m-1$,
\[\varphi_w(P_n, Q)=\sum_{k=0}^{m-1} \lambda_k \varphi_w(P_n, P_k)=0\]
donc $P_n\in\R_{m-1}[X]^{\perp}$. \hfill $\square$

\bigskip

\begin{rem} --- Fixons une fonction poids $w$ sur $I$ et munissons $\R[X]$ du produit scalaire $\varphi_w$. L'existence d'une famille de polynômes $w-$orthogonaux est assurée par le procédé de Gram-Schmidt qui permet de construire une telle famille à partir de la base canonique de $\R[X]$. On obtient ainsi une famille orthonormée de polynômes $(P_n)$. Si $(Q_n)$ est une autre famille de polynômes $w-$orthogonaux alors, pour tout $n\in\N^*$, $Q_n$ appartient au supplémentaire orthogonal de $\R_{n-1}[X]$ dans $\R_n[X]$ qui est une droite engendrée par $P_n$ donc il existe un réel $\alpha_n$ tel que $Q_n=\alpha_nP_n$. Ceci reste évidemment vraie pour $n=0$ puisque $P_0$ et $Q_0$ sont polynômes constants non nuls. Ainsi, pour tout $n\in\N$, le polynôme de degré $n$ d'une famille de polynômes $w-$orthogonaux est unique à une constante multiplicative près (constante dépendant de $n$). En particulier, il existe une unique famille de polynômes unitaires $w-$orthogonaux.
\label{rem_unicité_poly_ortho}
\end{rem}

\begin{propriete} --- Soit $w$ une fonction poids sur $I$ et $(P_n)$ une famille de polynômes $w-$orthogonaux. Alors, pour tout $n\in\N^*$, le polynôme $P_n$ est scindé à racines simples sur $\R$ et toutes ses racines appartiennent à $I$.
\label{prop_poly_ortho_racines_simples}
\end{propriete}

\emph{Preuve}. --- Soit $n\in\N^*$. Notons $r\in\llbracket 0, n\rrbracket$ le nombre de racines réelles distinctes de $P_n$ qui appartiennent à $I$ et qui sont de multiplicité impaire. Si $r=0$, posons $Q=1$ et, sinon, notons ces $r$ racines $x_1$, $x_2$, ..., $x_r$ et posons $Q=(X-x_1)(X-x_2)\cdots(X-x_r)$. Ainsi, $P_nQ$ est un polynôme dont toutes les racines réelles appartenant à $I$ sont de multiplicité paire. Par suite, $P_nQ$ est de signe constant sur $I$ et donc la fonction $t\mapsto P_n(t)Q(t)w(t)$ également. Comme celle-ci est continue et non identiquement nulle sur $I$, 
\[\varphi_w(P_n,Q) = \int_I P_n(t)Q(t)w(t) \diff t \neq 0\]
donc $P_n$ et $Q$ ne sont pas orthogonaux. Or, d'après le lemme \ref{lem_base_ortho}, $P_n\in\R_{n-1}[X]^{\perp}$ donc $Q\notin\R_{n-1}[X]$ et ainsi $\deg Q \geqslant n$. Or, $\deg Q=r \leqslant n$ donc $r=\deg Q=n$ et dès lors $P$ admet $n$ racines réelles distinctes dans $I$ et celles-ci sont donc toutes simples. \hfill $\square$

\bigskip

Dans la suite, nous nous intéresserons (entre autres) aux polynômes orthogonaux dits \og classiques \fg{} qui correspondent aux fonctions poids suivantes:
\begin{enumerate}
\item[$\bullet$] $w : t \mapsto (1-t)^{\alpha}(1+t)^{\beta}$ (avec $\alpha$ et $\beta$ deux réels strictement supérieurs à $-1$) définie sur $I=\intervalleoo{-1}{1}$ (polynômes de Jacobi);
\item[$\bullet$] $w : t \mapsto t^{\alpha}\mathrm{e}^{-t}$ (avec $\alpha$ un réel strictement supérieur à $-1$) définie sur $I=\intervalleoo{0}{+\infty}$ (polynômes de Laguerre);
\item[$\bullet$] $w : t \mapsto \mathrm{e}^{-t^2}$  définie sur $I=\R$ (polynômes d'Hermite).
\end{enumerate} 
Les familles de fonctions polynômiales associées à ces familles de polynômes apparaissent comme solutions d'équations différentielles dites hypergéométriques:
\[(E_n):~\sigma(x)y''+\tau(x)y'+\lambda_n y=0\]
où $n\in\N$, $\sigma$ et $\tau$ sont des fonctions polynômiales indépendantes de $n$ de degrés respectifs au plus 2 et 1 et $\lambda_n$ est un réel ne dépendant que de $n$. En 1929, Bochner \cite{Boc29} a déterminé les familles $(f_n)$ de fonctions polynomiales telles que, pour tout $n\in\N$, $f_n$ est une solution de $(E_n)$ de degré $n$. Celles-ci s'expriment à l'aide des polynômes de Jacobi, de Laguerre et d'Hermite ainsi qu'à l'aide d'un quatrième type de polynôme: les polynômes de Bessel. Ces derniers sont parfois également considérés comme une famille de polynômes orthogonaux mais dans un sens différent de celui que nous avons défini plus haut (voir \cite{KF49}). Sur toutes ces questions, on pourra consulter \cite[p. 105 et suiv.]{Pra10}, \cite[chap. 1]{NSU91} et \cite[chap. 1]{NU83}.

\subsection{Symbole de Pochhammer et coefficient binomial généralisé}
\begin{definition} --- Soit $z\in\C$ et $k\in\N$. On définit le symbole de Pochhammer, noté $(z)_k$, par 
\[(z)_k=\prod_{j=0}^{k-1} (z+j)\]
avec la convention habituelle $(z)_0=1$. On définit, de plus, le coefficient binomial généralisé $\binom{z}{k}$ par 
\[\binom{z}{k}=\dfrac{z(z-1)\cdots(z-k+1)}{k!}=\dfrac{(z-k+1)_k}{k!}.\]
\end{definition}

\begin{lemme} --- Soit $n\in\N$ et $k\in\llbracket 0, n\rrbracket$. Alors,
\[\binom{-\frac{1}{2}+n}{k}=\dfrac{(n-k)!(2n)!}{4^k k!n!(2n-2k)!} \qquad \text{et} \qquad \binom{\frac{1}{2}+n}{k}=\dfrac{(n-k+1)!(2n+2)!}{4^k k!(n+1)!(2n-2k+2)!}\]
\label{lem_coeff_binomial}
\end{lemme}

\emph{Preuve}. --- Par définition,
\begin{align*}
\binom{-\frac{1}{2}+n}{k}&=\dfrac{\left(-\frac{1}{2}+n\right)\left(-\frac{1}{2}+n-1\right)\cdots\left(-\frac{1}{2}+n-k+2\right)\left(-\frac{1}{2}+n-k+1\right)}{k!} \\
&=\dfrac{(2n-1)(2n-3)\cdots(2n-2k+3)(2n-2k+1)}{2^k k!} \\
&=\dfrac{(2n)!}{2^k k!(2n)(2n-2)\cdots(2n-2k+4)(2n-2k+2)(2n-2k)!} \\
&=\dfrac{(n-k)!(2n)!}{4^k k!n!(2n-2k)!}
\end{align*}
et, de plus,
\[\binom{\frac{1}{2}+n}{k}=\binom{-\frac{1}{2}+n+1}{k}=\dfrac{(n-k+1)!(2n+2)!}{4^k k!(n+1)!(2n-2k+2)!}.\]
\hfill $\square$

\subsection{Propriétés de la valuation $\boldsymbol{p-}$adique}
\label{para_valuation_p_adique}

\begin{lemme} --- Pour tout entier $d\geqslant 3$ et tout nombre premier $p$, $v_p(d)\leqslant d-2$.
\label{lem_majoration_vp}
\end{lemme}

\emph{Preuve}. --- Soit un entier $d\geqslant 3$ et un nombre premier $p$. Alors, comme $p\geqslant 2$,
\[p^{d-1}-1 =(p-1)\sum_{j=0}^{d-2} p^j \geqslant \sum_{j=0}^{d-2} p^j > \sum_{j=0}^{d-2} 1 =   d-1\]
donc $p^{d-1} > d$ ce qui assure que $v_p(d) < d-1$ i.e. $v_p(d) \leqslant d-2$. \hfill $\square$

\begin{lemme} --- Soit $s$ et $t$ deux entiers premiers entre eux tels que $t\geqslant 1$ et si $t=1$ alors $s\geqslant 0$. Pour tout nombre premier $p$ et tout entier $k\geqslant 1$,
\label{lem_majoration_factorielle}
\[v_p((s+t)(s+2t)\cdots(s+kt)) \leqslant v_p((\abs{s}+kt)!).\]
\end{lemme}

\emph{Preuve}. --- Soit $k\in\N^*$ et $P_k=(s+t)(s+2t)\cdots(s+kt)$. Il est équivalent de montrer que  $P_k$ divise $(\abs{s}+kt)!$. Notons que, pour tout $j\in\llbracket 1, k\rrbracket$, $s+jt\neq0$. En effet, dans le cas contraire, $\frac{s}{t}=-j$ donc, comme cette fraction est irréductible, $s=-j$ et $t=1$, ce qui est exclu puisque $s\geqslant 0$ si $t=1$. 

Si $(s+t)(s+kt) > 0$, $\abs{P_k}=\abs{s+t}\abs{s+2t}\cdots\abs{s+kt}$ est le produit de $k$ entiers distincts compris entre $1$ et $\max(\abs{s+t},\abs{s+kt}) \leq \abs{s}+kt$ car $k\geqslant 1$. Ainsi, dans ce cas, $P_k$ divise $(\abs{s}+kt)!$. Notons que cette condition est remplie si $t=1$ car alors $s\geqslant 0$. Dans la suite, on peut donc supposer $t\geqslant 2$.

Supposons à présent que $(s+t)(s+kt)<0$ i.e. que $s+t<0$ et $s+kt>0$ car $t>0$. S'il existe deux entiers $j$ et $\ell$ avec $1 \leqslant j < \ell\leqslant k$ tels que $\abs{s+j t}=\abs{s+\ell t}$ alors $s+\ell t=-s-jt$ donc $(\ell+j)t=-2s$. Comme $t$ est premier avec $s$, ceci implique que $t$ divise $2$ donc $t=2$. Ainsi, si $t\geqslant3$, comme précédemment, $\abs{P_k}$ est le produit de $k$ entiers distincts compris entre $1$ et $\abs{s}+kt$ donc $P_k$ divise $(\abs{s}+kt)!$. Reste le cas où $t=2$. Dans ce cas, comme $s$ est premier avec $t$ et $s<0$ (car $s+t<0$), il existe un entier $m\geqslant 0$ tel que $s=-2m-1$. Ainsi, 
\[P_k=(-2m+1)\times(-2m+3)\times \cdots\times (-3) \times (-1) \times 1 \times 3 \times \cdots \times (2(k-m)-1)\]
donc, puisque $(-2m+1)\times(-2m+3)\times \cdots \times (-3) \times (-1)$ divise $2\times 4 \times 6 \times \cdots \times (4m-2)$, en notant $M=\max(4m-2,2(k-m)-1)$, $P_k$ divise $M!$. Or, d'une part, $s+kt>0$ donc $kt>-s=\abs{s}$ et ainsi $\abs{s}+kt>2\abs{s}=4m+2>4m-2$ et, d'autre part, $\abs{s}+kt\geqslant 2k > 2(k-m)-1$. Il s'ensuit que $M \leqslant \abs{s}+kt$ donc $P_k$ divise $(\abs{s}+kt)!$. \hfill $\square$

\bigskip

Un des outils que nous utiliserons souvent dans la suite est la formule suivante, bien connue, due à Legendre \cite[XVI p. 8]{Leg08}.

\begin{propriete}[Legendre, 1808] --- Pour tout entier $n\geqslant 1$ et tout nombre premier $p$, la valuation $p-$adique de $n!$ est 
\[v_p(n!)=\sum_{j=1}^{+\infty} \entiere{\frac{n}{p^j}}=\sum_{j=1}^{r_n} \entiere{\dfrac{n}{p^j}}\]
où $r_n=\entiere{\frac{\ln n}{\ln p}}$.
\label{prop_formule_Legendre}
\end{propriete}

\emph{Preuve}. --- Soit un entier $n\geqslant 1$ et $p$ un nombre premier. Soit $j\in\N^*$. Les entiers $k\in\llbracket 1, n\rrbracket$ tels que $p^j$ divise $k$ sont $p^j$, $2p^j$, ..., $\entiere{\frac{n}{p^j}} p^j$: il y en a donc $\entiere{\frac{n}{p^j}}$. Le nombre d'entiers $k\in\llbracket 1, n\rrbracket$ tels que $v_p(k)=j$ i.e. le nombre d'entiers divisibles par $p^j$ mais pas par $p^{j+1}$ est donc $\entiere{\frac{n}{p^j}}-\entiere{\frac{n}{p^{j+1}}}$. On en déduit que
\[v_p(n!)=\sum_{k=1}^n v_p(k) = \sum_{j=1}^{+\infty} j\left(\entiere{\dfrac{n}{p^j}}-\entiere{\dfrac{n}{p^{j+1}}}\right).\]
Notons que cette somme est en fait finie puisque si $p^j >n$ i.e. si $j>\entiere{\frac{\ln n}{\ln p}}$ alors $\entiere{\frac{n}{p^j}}=0$. Dès lors,
\[v_p(n!)=\sum_{j=1}^{+\infty} j\entiere{\dfrac{n}{p^j}}-\sum_{j=1}^{+\infty} j\entiere{\dfrac{n}{p^{j+1}}} = \sum_{j=1}^{+\infty} j\entiere{\dfrac{n}{p^j}}-\sum_{j=2}^{+\infty} (j-1)\entiere{\dfrac{n}{p^{j}}}=\entiere{\dfrac{n}{p}} + \sum_{j=2}^{+\infty} [j-(j-1)]\entiere{\dfrac{n}{p^j}}\]
et donc
\[v_p(n!)=\sum_{j=1}^{+\infty} \entiere{\dfrac{n}{p^j}}=\sum_{j=1}^{r_n} \entiere{\dfrac{n}{p^j}}.\]
 \hfill $\square$

\begin{corollaire} --- Soit un entier naturel $n\geqslant 1$. Pour tout nombre premier $p$, $v_p(n!) < \frac{n}{p-1}$.
\label{coro_Legendre}
\end{corollaire}

\emph{Preuve}. --- Soit un nombre premier $p$. En posant $r_n=\entiere{\frac{\ln n}{\ln p}}$, d'après la formule de Legendre,
\[v_p(n!)=\sum_{j=1}^{r_n} \entiere{\dfrac{n}{p^j}} \leqslant \sum_{j=1}^{r_n} \dfrac{n}{p^j} <  n\sum_{j=1}^{+\infty} \dfrac{1}{p^j} = \dfrac{n}{p-1}.\]
\hfill $\square$ 

\subsection{Critères d'irrationalité}
\label{para_critère_irrationalité}

\begin{propriete} --- Soit $P=\sum\limits_{k=0}^d c_kX^k \in\Z[X]$ un polynôme de degré $d>0$. Si un nombre rationnel $\alpha$ écrit sous forme irréductible $\alpha=\frac{s}{t}$ est racine de $P$ alors $s$ divise $c_0$ et $t$ divise $c_d$.
\label{prop_RRT}
\end{propriete}

\emph{Preuve}. --- Supposons qu'un rationnel $\alpha=\frac{s}{t}$ écrit sous forme irréductible soit racine de $P$. Comme $\sum\limits_{k=0}^d c_k \left(\frac{s}{t}\right)^k=0$, en multipliant par $t^d$, on obtient $\sum\limits_{k=0}^d c_k s^kt^{d-k}=0$. En particulier, on a 
\[c_ds^d=-\sum\limits_{k=0}^{d-1} c_k s^kt^{d-k}=-t\sum\limits_{k=0}^{d-1} c_k s^kt^{d-k-1}.\]
\'Etant donné que $-\sum\limits_{k=0}^{d-1} c_k s^kt^{d-k-1}$ est un entier, il s'ensuit que $t$ divise $c_ds^d$. Or, comme $s$ et $t$ sont premiers entre eux, il en est de même de $s^d$ et $t$ donc, d'après le lemme de Gauss, $t$ divise $c_d$. De même, on a $c_0t^d=-\sum\limits_{k=1}^{d} c_k s^{k}t^{d-k}=-s\sum\limits_{k=1}^{d} c_k s^{k-1}t^{d-k}$ donc $s$ divise $c_0$. \hfill $\square$

\begin{definition} --- On dit qu'un nombre complexe $\alpha$ est un \emph{entier algébrique} s'il existe un polynôme unitaire $P\in\Z[X]$ tel que $P(\alpha)=0$.
\end{definition}

On déduit immédiatement de la proposition \ref{prop_RRT} le corollaire suivant.

\begin{corollaire} --- Si un entier algébrique est rationnel alors il est entier.
\label{coro_entier_alg}
\end{corollaire}

\begin{propriete} --- Soit un entier $d\geqslant 2$. On considère des entiers $b_0$, $b_1$, ..., $b_d$ tels que $b_0\neq 0$. On suppose qu'il existe un nombre premier $p$ tel que $p$ divise $b_k$ pour tout $k\in\llbracket 0, d-1\rrbracket$ et tel que $p$ ne divise pas $b_d$. On suppose, de plus, que, pour tout $k\in\llbracket 1, d\rrbracket$, $v_p(b_k) > v_p(b_0)-k$. Alors, pour tous entiers $a_0$, $a_1$, ..., $a_d$ tels que $p$ ne divise ni $a_0$ ni $a_d$, le polynôme $P:=\sum\limits_{k=0}^{d} a_kb_kX^k$ n'admet pas de racine rationnelle. 
\label{prop_cri_irra_poly_Newton}
\end{propriete}

\emph{Preuve}. --- Soit $a_0$, $a_1$, ..., $a_d$ des entiers tels que $p$ ne divise ni $a_0$ ni $a_d$. Posons $P:=\sum\limits_{k=0}^{d} a_kb_kX^k$ et supposons que $\alpha=\frac{s}{t}$ soit une racine rationnelle de $P$ écrite sous forme irréductible. Alors, $ 0=t^dP(\alpha)=\sum\limits_{k=0}^{d} a_kb_ks^kt^{d-k}$ donc $a_db_ds^d=-\sum\limits_{k=0}^{d-1} a_kb_ks^kt^{d-k}$. Or, pour tout $k\in\llbracket 0, d-1\rrbracket$, $p$ divise $b_k$ donc $p$ divise $a_db_ds^d$. Mais, comme, par hypothèse, $p$ ne divise ni $a_d$ ni $b_d$, on conclut que $p$ divise $s$. 

Notons $m=v_p(b_0)$. Alors, comme précédemment, $P(\alpha)=0$ assure que $a_0b_0t^d=-\sum\limits_{k=1}^{d} a_kb_ks^kt^{d-k}$ donc
\begin{equation}
a_0\dfrac{b_0}{p^m}t^{d}=-\sum_{k=1}^{d} a_k\dfrac{b_ks^k}{p^m}t^{d-k}.
\label{eq_prop_irra}
\end{equation}
Or, d'une part, comme $s$ et $t$ sont premiers entre eux, $p$ ne divise pas $t$ et, par hypothèse, $p$ ne divise pas $a_0$ donc, par définition de $m$, $a_0\frac{b_0}{p^m}t^{d}$ est un entier premier avec $p$. D'autre part, pour tout $k\in\llbracket 1, d\rrbracket$,
\[v_p(b_ks^k) =v_p(b_k)+kv_p(s) > m-k+kv_p(s) \geqslant m\]
car $v_p(s)\geqslant 1$. Ainsi, pour tout $k\in\llbracket 1, d\rrbracket$, $\frac{b_ks^k}{p^m}$ est un entier divisible par $p$ donc $\sum\limits_{k=1}^{d} a_k\frac{b_ks^k}{p^m}t^{d-k}$ est également un entier divisible par $p$. L'égalité \eqref{eq_prop_irra} conduit alors à une contradiction donc $P$ n'a pas de racine rationnelle.  \hfill $\square$

\bigskip

\begin{rem} --- De manière plus générale, on peut démontrer le résultat suivant:

\medskip

\emph{Soit $Q=\sum\limits_{k=0}^{d} b_kX^k\in\Z[X]$ un polynôme de degré $d \geqslant 2$ tel que $b_0\neq 0$. On suppose qu'il existe un entier $\ell\in\llbracket 0, d-1 \rrbracket$ et un nombre premier $p$ tel que $p$ divise $b_k$ pour tout $k\in\llbracket 0, d-\ell-1\rrbracket$ et $p$ ne divise pas $b_d$. On suppose, de plus, qu'il existe un entier $h>\ell$ tel que, pour tout $k\in\llbracket 1, d\rrbracket$, $\frac{v_p(b_k)-v_p(b_0)}{k}>-\frac{1}{h}$. Alors, pour tous entiers $a_0$, $a_1$, ..., $a_d$ tels que $p$ ne divise ni $a_0$ ni $a_d$, le polynôme $P:=\sum\limits_{k=0}^{d} a_kb_kX^k$ n'admet pas de diviseur dans $\Z[X]$ de degré $s\in\llbracket \ell+1, h\rrbracket$.}

\medskip

La propriété \ref{prop_cri_irra_poly_Newton} correspond au cas $\ell=0$ et $h=1$. Ce théorème général est utilisé dans de nombreuses démonstrations d'irréductibilité (voir par exemple \cite{Fil95}, \cite{Fil96} et \cite{FT02}). Le nombre $\min\limits_{k\in\llbracket 1, d\rrbracket} \frac{v_p(b_k)-v_p(b_0)}{k}$ s'interprète comme la pente du premier segment du polygone de Newton de $Q$ relativement à $p$ et la démonstration du théorème s'appuie sur le théorème de Dumas. Sur ces questions, on pourra consulter \cite[p. 52 et suiv.]{Pra10}.
\label{rem_poly_Newton}
\end{rem}

\section{Deux exemples de polynômes de Jacobi}

\subsection{Polynômes de Tchebichef}
\label{secTchebichef}
\subsubsection{Définition et propriétés}

\begin{definition} --- On considère les suites $(T_n)$ et $(U_n)$ d'éléments de $\Z[X]$ définies par
\[\begin{cases} T_0=1,~T_1=X, \\ \forall n\in\N,~T_{n+2}=2XT_{n+1}-T_n \end{cases} \qquad\text{et} \qquad \begin{cases} U_0=1,~U_1=2X, \\ \forall n\in\N,~U_{n+2}=2XU_{n+1}-U_n \end{cases}\]

Pour tout $d\in\N$, $T_d$ est appelé le polynôme de Tchebichef de première espèce d'indice $d$ et $U_d$ le polynôme de Tchebichef de seconde espèce d'indice $d$. 
\end{definition}

Par une récurrence double immédiate, on montre que, pour tout $d\in\N$, $\deg T_d=\deg U_d=d$, que le coefficient dominant de $T_d$ est $2^{d-1}$ pour tout $d\in\N^*$ et que celui de $U_d$ est $2^d$ pour tout $d\in\N$.

\begin{lemme} --- Soit $x\in\R$. Pour tout $d\in\N$,
\[T_d(\cos(x))=\cos(dx) \qquad \text{et} \qquad \sin(x)U_d(\cos(x))=\sin((d+1)x).\]
\label{lem_Tchebichef}
\end{lemme}

\emph{Preuve}. --- On raisonne par récurrence double. Le résultat est immédiat pour $d=0$ et $d=1$ (puisque $2\sin(x)\cos(x)=\sin(2x)$). Supposons les égalités établies aux rangs $d$ et $d+1$. Alors,
\begin{align*}
T_{d+2}(\cos (x))&=2\cos(x)T_{d+1}(\cos(x))-T_d(\cos(x))=2\cos(x)\cos((d+1)x)-\cos(dx) \\
&=\cos((d+1)x-x)+\cos((d+1)x+x)-\cos(dx)=\cos((d+2)x
\end{align*}
et
\begin{align*}
\sin(x)U_{d+2}(\cos (x))&=\sin(x)\left[2\cos(x)U_{d+1}(\cos(x))-U_d(\cos(x))\right]=2\cos(x)\sin((d+2)x)-\sin((d+1)x) \\
&=\sin((d+2)x-x)+\sin((d+2)x+x)-\sin((d+1)x)=\sin((d+3)x
\end{align*}

Ainsi, les égalités sont établies aux rangs $d+1$ et $d+2$, ce qui achève la récurrence double. \hfill $\square$

\bigskip

\begin{rem} ---  Il est simple de voir que $w_1 : t \mapsto \frac{1}{\sqrt{1-t^2}}$ et $w_2 : t \mapsto \sqrt{1-t^2}$ définissent des fonctions poids sur $\intervalleoo{-1}{1}$ et que $(T_n)$ est une famille de polynômes $w_1-$orthogonaux et $(U_n)$ une famille de polynômes $w_2-$orthogonaux. Ainsi, d'après la propriété \ref{prop_poly_ortho_racines_simples}, pour tout $n\in\N^*$, $T_n$ et $U_n$ sont scindés à racines simples sur $\R$ et toutes leurs racines appartiennent à $\intervalleoo{-1}{1}$. Cependant, dans le cas des polynômes de Tchebichef, il est encore plus simple d'expliciter ces racines.
\label{rem_ortho_Tchebichef}
\end{rem}

\begin{propriete} --- Soit $d\in\N^*$. Les racines de $T_d$ sont les réels $a_k:=\cos\left(\frac{2k-1}{2d}\pi\right)$ pour $k\in\llbracket 1, d\rrbracket$ et les racines de $U_d$ sont les réels $b_k:=\cos\left(\frac{k}{d+1}\pi\right)$ pour $k\in\llbracket 1, d\rrbracket$.
\end{propriete}

\emph{Preuve}. --- Pour tous entiers $k$ et $\ell$ tels que $1 \leqslant k < \ell \leqslant d$, $0 < \frac{2k-1}{2d} \pi< \frac{2\ell-1}{2d} \pi < \pi$ et $0  < \frac{k}{d+1}  \pi< \frac{\ell}{d+1} \pi < \pi$ donc, par injectivité de cos sur $\intervalleoo{0}{\pi}$, $a_k \neq a_{\ell}$ et $b_k \neq b_{\ell}$. Ainsi, les $d$ réels $a_k$ sont tous distincts et les $d$ réels $b_k$ sont tous distincts. 

Soit $k\in\llbracket 1, d\rrbracket$. Alors,
\[T_d\left(a_k\right)=\cos\left(d\dfrac{2k-1}{2d}\pi\right)=\cos\left(k\pi-\dfrac{\pi}{2}\right)=0\]
donc $a_k$ est bien racine de $T_d$.
Par ailleurs, 
\[\sin\left(\dfrac{k}{d+1}\pi\right)U_d\left(b_k\right)=\sin\left((d+1)\dfrac{k}{d+1}\pi\right)=\sin\left(k\pi\right)=0\]
et, comme $0<\frac{k}{d+1}\pi<\pi$, $\sin\left(\frac{k}{d+1}\pi\right)\neq 0$ donc $U_d\left(b_k\right)=0$.

Comme $T_d$ et $U_d$ sont de degré $d$, la conclusion s'ensuit.  \hfill $\square$

\subsubsection{Irrationalité des racines}

\'Etudier l'irrationalité des racines de polynômes de Tchebichef se ramène donc à étudier l'irrationalité de $\cos(r\pi)$ pour $r\in\Q$. Ce sujet est abondamment traité dans la littérature (voir, par exemple, les références données dans \cite[p. 41]{Ni56}). Nous suivons ici la démonstration de Maier \cite{Mai65} (voir également \cite[p. 308]{NZM91} et \cite[p. 103]{Pra10})

\begin{propriete} --- Pour tout $r\in\Q$, $\cos(r\pi)$ est rationnel si et seulement si $r \in \frac{1}{2}\Z\cup\frac{1}{3}\Z$.
\label{prop_irrationalité_cos}
\end{propriete}

\emph{Preuve}. --- Considérons la suite de polynômes $(P_n)_{n\in\N^*}$ définie, pour tout $n\in\N^*$, par $P_n=2T_n\left(\frac{X}{2}\right)$. Alors, $P_1=X$, $P_2=X^2-2$ et, pour tout $n\in\N^*$,
\[P_{n+2}=2T_{n+2}\left(\dfrac{X}{2}\right)=4\left(\dfrac{X}{2}\right)T_{n+1}\left(\dfrac{X}{2}\right)-2T_n\left(\dfrac{X}{2}\right)=XP_{n+1}-P_n\]
donc, par une récurrence double immédiate, pour tout $n\in\N^*$, $P_n$ est un polynôme unitaire de $\Z[X]$. 

Soit $r=\frac{s}{t}$ un rationnel écrit sous forme irréductible. Alors, d'après le lemme \ref{lem_Tchebichef},
\[P_t(2\cos(r\pi))=2T_t(\cos(r\pi))=2\cos(tr\pi)=2\cos(s\pi).\]
Ainsi, $2\cos(r\pi)$ est racine du polynôme $F:=P_t-2\cos(s\pi)$. Comme $\cos(s\pi)\in\Z$, $F$ est un polynôme unitaire de degré $t\geqslant 1$ de $\Z[X]$ donc $2\cos(r\pi)$ est un entier algébrique.

Supposons que $\cos(r\pi)$ est rationnel. Alors, il en est de même de $2\cos(r\pi)$ et donc, d'après la propriété \ref{prop_RRT}, $2\cos(r\pi)$ est entier. De plus, comme $-2 \leqslant 2\cos(r\pi) \leqslant 2$, $\cos(r\pi)\in\left\{-1, -\frac{1}{2}, 0, \frac{1}{2}, 1\right\}$.

Si $\cos(r\pi)\in\{-1, 0, 1\}$ alors $r\pi\equiv 0 \pmod{\frac{\pi}{2}}$ donc $r\in\frac{1}{2}\Z$.

Sinon, $\cos(r\pi)\in\left\{-\frac{1}{2}, \frac{1}{2}\right\}$ donc $r\pi\equiv \frac{\pi}{3} \pmod{\pi}$ ou $r\pi\equiv \frac{2\pi}{3} \pmod{\pi}$. Par suite, il existe un entier $k$ tel que $r=\frac{1}{3}+k=\frac{1+3k}{3}$ ou $r=\frac{2}{3}+k=\frac{2+3k}{3}$ et donc $r\in\frac{1}{3}\Z$. 

On a donc démontré que si $r$ est un rationnel tel que $\cos(r\pi)$ est également un nombre rationnel alors $r\in\frac{1}{2}\Z\cup\frac{1}{3}\Z$. La réciproque est immédiate. \hfill $\square$

\begin{theoreme} --- Soit un entier $d\geqslant 1$. 
\begin{enumerate}
\item Si $d$ est pair, les racines de $T_d$ sont toutes irrationnelles et, si $d$ est impair, l'unique racine rationnelle de $T_d$ est $0$.
\item Les seules racines rationnelles possibles pour $U_d$ sont $0$, $-\frac{1}{2}$ et $\frac{1}{2}$. De plus, $0$ est racine de $U_d$ si et seulement si $d$ est impair et $\frac{1}{2}$ et $-\frac{1}{2}$ sont racines de $U_d$ si et seulement si $d\equiv 2 \pmod{3}$.
\end{enumerate} 
\end{theoreme}

\emph{Preuve}. --- Soit $k\in\llbracket 1, d\rrbracket$.

\begin{enumerate}
\item Notons $r_k:=\frac{2k-1}{2d}$ et $a_k=\cos\left(r_k\pi\right)$. La propriété \ref{prop_irrationalité_cos} assure que $a_k$ est rationnel si et seulement si $2r_k$ ou $3r_k$ est un entier. Comme $3(2k-1)$ est impair, $2d$ ne divise pas $3(2k-1)$ donc $3r_k \notin \Z$. De plus, comme $0 < 2k-1 < 2d$, $2r_k=\frac{2k-1}{d}$ est entier si et seulement si $2k-1=d$ i.e. $k=\frac{d+1}{2}$. Dans ce cas, qui est possible si et seulement si $d$ est impair, $r_k=\frac{1}{2}$ et donc $a_k=0$. 

\item De même, notons $t_k:=\frac{k}{d+1}$ et $b_k=\cos\left(t_k\pi\right)$. La propriété \ref{prop_irrationalité_cos} assure que $b_k$ est rationnel si et seulement si $2t_k$ ou $3t_k$ est entier. Or, comme $0 < k < d$, $0< 2t_k < 2$ donc $2t_k$ est entier si et seulement si $2t_k=1$ i.e. $t_k=\frac{1}{2}$. Dans ce cas, qui se produit si et seulement si $d$ est impair et $k=\frac{d+1}{2}$, $a_k=0$. De même, $0 < 3t_k < 3$ donc  $3t_k$ est entier si et seulement si $t_k=\frac{1}{3}$ ou $t_k=\frac{2}{3}$. Ces deux cas se produisent si et seulement si $d\equiv 2 \pmod{3}$ avec $k=\frac{d+1}{3}$ ou $k=\frac{2d+2}{3}$ et on trouve alors $a_{k}=\frac{1}{2}$ ou $a_{k}=-\frac{1}{2}$. 
\end{enumerate}
\hfill $\square$

\bigskip

\begin{rem} --- On peut montrer que $T_d$ est irréductible sur $\Q$ si et seulement si $d$ est une puissance de $2$ et $U_d$ est irréducible sur $\Q$ si et seulement si $d=1$. De plus, si $d$ est impair alors $\frac{T_ d}{X}$ est irréductible sur $\Z$ si et seulement si $d$ est premier. À ce sujet, on pourra consulter \cite{RTW05} qui donne une description complète de la factorisation de $T_d$ et $U_d$ en produit d'irréductibles dans $\Z[X]$.
\label{rem_Tchebichef_irré}
\end{rem}

\subsection{Polynômes de Legendre}

\subsubsection{Définition et propriétés}

\begin{definition} --- Pour tout $d\in\N$, On pose $\mathcal{L}_d=(X^2-1)^d$ et on définit le polynôme de Legendre d'indice $d$ par
\[P_d:= \dfrac{1}{2^d d!} \mathcal{L}_d^{(d)}.\]
\end{definition}

Il est clair que, pour tout $d\in\N$, $P_d$ est un polynôme de $\Q[X]$ de degré $d$.

\begin{propriete} --- Soit $w$ la fonction constante égale à $1$ sur l'intervalle $I:=\intervalleoo{-1}{1}$. Alors, la famille $(P_n)_{n\in\N}$ est une famille de polynômes $w-$orthogonaux.
\end{propriete}

\emph{Preuve}. --- Comme on l'a déjà dit, pour tout $n\in\N$, $P_n$ est de degré $n$.

Soit $m$ et $n$ deux entiers naturels tels que $m<n$. Alors, par intégrations par parties itérées,
\begin{align*}
\varphi_{w}(P_m,P_n)&=\int_{-1}^1 P_m(t)P_n(t) \diff t=\dfrac{1}{2^n n!}\int_{-1}^1 P_m(t)\mathcal{L}_n^{(n)}(t) \diff t \\
&=\dfrac{1}{2^n n!}\sum_{k=0}^{n-1} (-1)^k\left[P_m^{(k)}(t)\mathcal{L}_n^{(n-k-1)}(t)\right]_{-1}^1+\dfrac{(-1)^n}{2^n n!} \int_{-1}^{1} P_m^{(n)}(t)\mathcal{L}_n(t) \diff t
\end{align*}
Comme $P_m$ est de degré $m<n$, $P_m^{(n)}=0$ donc 
\[\int_{-1}^{1} P_m^{(n)}(t)\mathcal{L}_n(t) \diff t=0.\]
Par ailleurs, par définition, $1$ et $-1$ sont racines de $\mathcal{L}_n$ de multiplicité $n$ donc ils sont également racines de $\mathcal{L}_n^{(n-k-1)}$ pour tout $k\in\llbracket 0, n-1 \rrbracket$ et ainsi, pour tout $k\in\llbracket 0, n-1 \rrbracket$, $\left[P_m^{(k)}(t)\mathcal{L}_n^{(n-k-1)}(t)\right]_{-1}^1=0$. Il s'ensuit que $\varphi_{w}(P_m,P_n)=0$ et on conclut que $(P_n)_{n\in\N}$ est bien une famille de polynôme $w-$orthogonaux. \hfill $\square$

\bigskip

On déduit alors de la propriété \ref{prop_poly_ortho_racines_simples} le résultat suivant.

\begin{corollaire} --- Pour tout $d\in\N^*$, $P_d$ est scindé à racines simples sur $\R$ et toutes ses racines appartiennent à $\intervalleoo{-1}{1}$.
\label{coro_Legendre_scindé}
\end{corollaire}

\begin{propriete} --- Pour tout $d\in\N$, 
\[(1)~P_d=\dfrac{1}{2^d}\sum_{k=0}^d \dbinom{d}{k}^2(X+1)^k(X-1)^{d-k} \qquad \text{et} \qquad (2)~P_d=\dfrac{1}{2^d}\sum_{k=0}^{\entiere{\frac{d}{2}}} (-1)^k \dbinom{d}{k}\dbinom{2d-2k}{d} X^{d-2k}.\]
\label{prop_expression_Legendre}
\end{propriete}

\emph{Preuve}. --- Soit $d\in\N$. 

(1) En écrivant $\mathcal{L}_d=(X-1)^d(X+1)^d$, il vient, grâce à la formule de Leibniz,
\begin{align*}
P_d&=\dfrac{1}{2^d d!}\sum_{k=0}^{d} \dbinom{d}{k} \left[(X-1)^d\right]^{(k)}\left[(X+1)^d\right]^{(d-k)} \\
&=\dfrac{1}{2^d d!}\sum_{k=0}^{d} \dbinom{d}{k} \dfrac{d!}{(d-k)!}(X-1)^{d-k}\dfrac{d!}{k!}(X+1)^{k} \\
&=\dfrac{1}{2^d}\sum_{k=0}^{d} \dbinom{d}{k} \dfrac{d!}{k!(d-k)!}(X+1)^{k}(X-1)^{d-k}
\end{align*}
donc $P_d=\frac{1}{2^d}\sum_{k=0}^d \binom{d}{k}^2(X+1)^k(X-1)^{d-k}$.

(2) D'après la formule du binôme de Newton, $\mathcal{L}_d=\sum\limits_{k=0}^d \binom{d}{k}(-1)^kX^{2d-2k}$. Notons, pour tout $k\in\llbracket 0, d \rrbracket$, $M_k:=X^{2d-2k}$. Alors, si $2d-2k<d$ i.e. si $k>\entiere{\frac{d}{2}}$, $M_k^{(d)}=0$ et, sinon, 
\[M_k^{(d)}=\frac{(2d-2k)!}{(d-2k)!}X^{d-2k}=d!\binom{2d-2k}{d}X^{d-2k}.\]
Ainsi,
\[P_d=\dfrac{1}{2^d d!}\sum_{k=0}^{\entiere{\frac{d}{2}}} \binom{d}{k}(-1)^kd!\binom{2d-2k}{d}X^{d-2k} = \dfrac{1}{2^d}\sum_{k=0}^{\entiere{\frac{d}{2}}} (-1)^k \binom{d}{k}\binom{2d-2k}{d}X^{d-2k}.\]
\hfill $\square$

\subsubsection{Irrationalité des racines}

\begin{theoreme} --- Soit un entier $d\geqslant 1$. Si $d$ est pair, les $d$ racines de $P_d$ sont irrationnelles et, si $d$ est impair, les $d-1$ racines non nulles de $P_d$ sont irrationnelles.
\end{theoreme}

\noindent \textbf{\emph{Démostration}} --- Comme, d'après le corollaire \ref{coro_Legendre_scindé}, $P_d$ est scindé sur $\R$, il suffit de montrer que la seule racine rationnelle possible de $P_d$ est $0$ et que, de plus, $0$ est racine de $P_d$ si et seulement si $d$ est impair.

Posons $Q_d:=(X-1)^dP_d\left(\frac{X+1}{X-1}\right)$. Alors, grâce à l'identité (1) de la propriété \ref{prop_expression_Legendre},
\begin{align*}
Q_d&=\dfrac{(X-1)^d}{2^d}\sum_{k=0}^d \dbinom{d}{k}^2 \left(\dfrac{X+1}{X-1}+1\right)^k\left(\dfrac{X+1}{X-1}-1\right)^{d-k} \\
&=\dfrac{(X-1)^d}{2^d}\sum_{k=0}^d \dbinom{d}{k}^2 \left(\dfrac{2X}{X-1}\right)^k\left(\dfrac{2}{X-1}\right)^{d-k} \\
&=\dfrac{(X-1)^d}{2^d}\sum_{k=0}^d \dbinom{d}{k}^2 \dfrac{2^d}{(X-1)^d}X^k
\end{align*}
soit $Q_d=\sum\limits_{k=0}^d \binom{d}{k}^2 X^k$. En particulier, $Q_d$ est un polynôme de $\Z[X]$.

Supposons que $r$ soit une racine rationnelle de $P_d$. D'après l'identité (1) de la propriété \ref{prop_expression_Legendre}, $P_d(1)=\frac{1}{2^d}\times (1+1)^d=1\neq 0$ donc $r\neq 1$. Posons alors $\alpha=\frac{r+1}{r-1}$ de sorte que $r=\frac{\alpha+1}{\alpha-1}$. Comme $r$ est rationnel, il en est de même pour $\alpha$. De plus, $Q_d(\alpha)=(\alpha-1)^dP_d(r)=0$ donc $\alpha$ est une racine rationnelle de $Q_d$. Or, $Q_d$ est un polynôme unitaire de $\Z[X]$ dont le coefficient constant est 1 donc, d'après la propriété \ref{prop_RRT}, $\alpha\in\{-1, 1\}$. Or, $\alpha\neq 1$ par définition donc $\alpha=-1$ i.e. $r=0$. Ainsi, la seule racine rationnelle possible de $P_d$ est $0$. De plus, d'après l'identité (2) de la propriété \ref{prop_expression_Legendre}, si $d$ est pair, en posant $\ell=\frac{d}{2}$, on a $P_{d}(0)=\frac{1}{2^d}(-1)^{\ell}\binom{d}{\ell} \neq 0$ donc $0$ n'est pas racine de $P_d$ dans ce cas. En revanche, si $d$ est impair alors la valuation de $P_d$ est $d-2\entiere{\frac{d}{2}}=d-(d-1)=1$ donc $0$ est racine de $P_d$. \hfill $\square$

\begin{rem} --- On conjecture que $P_d$ est irréductible sur $\Q$ pour tout entier pair $d$ et que $\frac{P_d}{X}$ est irréductible sur $\Q$ pour tout entier impair $d$. Ceci a été démontré dans certains cas particuliers (voir, par exemple, \cite{Hol12}, \cite{Hol13}, \cite{Wah52} et \cite{Wah60}) mais le cas général reste un problème ouvert.
\label{rem_Legendre_irré}
\end{rem}

\subsection{Quelques mots sur le cas général}

\begin{definition} --- Pour tous réel $\alpha>-1$ et $\beta>-1$ et tout entier naturel $d$, on définit le polynôme de Jacobi d'indice $d$ et de paramètres $\alpha$ et $\beta$ par 
\[P_{d}^{(\alpha,\beta)}:=\dfrac{1}{2^d}\sum_{k=0}^{d} \binom{\alpha+d}{d-k}\binom{\beta+d}{k}(X-1)^{k}(X+1)^{d-k}.\]
\end{definition}

Il est clair que $P_d^{(\alpha,\beta)}$ est un polynôme de degré $d$ de $\R[X]$ (dont le coefficient dominant est $\frac{1}{2^d}\sum\limits_{k=0}^{d} \binom{\alpha+d}{d-k}\binom{\beta+d}{k}>0$) et que, de plus, si $\alpha$ et $\beta$ sont rationnels alors $P_d^{(\alpha,\beta)}$ appartient à $\Q[X]$. 

On vérifie facilement, à l'aide de la propriété \ref{prop_expression_Legendre}, que, pour tout $d\in\N$, $P_{d}^{(0,0)}=P_d$ est le polynôme de Legendre d'indice $d$. En remarquant que, pour tous réels $\alpha$ et $\beta$ strictement supérieurs à $-1$ et tout entier naturel $d$, la dérivée d'ordre $d$ de la fonction $f_{d,\alpha,\beta}: t \mapsto (1-t)^{d+\alpha}(1+t)^{d+\beta}$ définie sur $\intervalleoo{-1}{1}$ est la fonction $f_{d,\alpha,\beta}^{(d)} : t\mapsto (-1)^d2^dd!(1-t)^{\alpha}(1+t)^{\beta}P_d^{(\alpha,\beta)}(t)$, on montre sans difficulté la propriété suivante.

\begin{propriete} --- Soit deux réels $\alpha$ et $\beta$ strictement supérieurs à $-1$. On considère la fonction $w_{\alpha, \beta} : t\mapsto (1-t)^{\alpha}(1+t)^{\beta}$ définie sur l'intervalle $I:=\intervalleoo{-1}{1}$. Alors, la famille $(P_n^{(\alpha, \beta)})_{n\in\N}$ est une famille de polynômes $w_{\alpha, \beta}-$orthogonaux.
\end{propriete}

On déduit alors de la propriété \ref{prop_poly_ortho_racines_simples} le résultat suivant.

\begin{corollaire} --- Pour tout $d\in\N^*$ et pour tous réels $\alpha>-1$ et $\beta>-1$, le polynôme $P_d^{(\alpha, \beta)}$ est scindé à racines simples sur $\R$ et toutes ses racines appartiennent à $\intervalleoo{-1}{1}$.
\label{coro_Jacobi_scindé}
\end{corollaire}

\begin{rem} --- On a $w_{-\frac{1}{2},-\frac{1}{2}}: t\mapsto \frac{1}{\sqrt{1-t^2}}$ donc, d'après les remarques \ref{rem_unicité_poly_ortho} et \ref{rem_ortho_Tchebichef}, pour tout $d\in\N$, il existe une constante $\lambda_d$ telle que $P_d^{(-\frac{1}{2},-\frac{1}{2})}=\lambda_d T_d$. En égalant les coefficients dominants, on déduit du lemme \ref{lem_coeff_binomial} que, pour tout $d\in\N$, $\lambda_d=\frac{(2d)!}{4^d(d!)^2}=\binom{-\frac{1}{2}+d}{d}$ donc $P_d^{(-\frac{1}{2},-\frac{1}{2})}=\binom{-\frac{1}{2}+d}{d} T_d$.

De la même façon, $w_{\frac{1}{2},\frac{1}{2}}: t\mapsto \sqrt{1-t^2}$ et on en déduit que, pour tout $d\in\N$, $P_d^{(\frac{1}{2},\frac{1}{2})}=2\binom{\frac{1}{2}+d}{d+1} U_d$
\label{rem_Tcebichef_Jacobi}
\end{rem}

L'étude des deux cas particuliers précédents montre qu'il existe des polynômes de Jacobi de degré arbitrairement grand ayant des racines rationnelles (cas des polynômes de Tchebichef) mais aussi des polynômes de Jacobi de n'importe quel degré (au moins égal à 2) n'ayant aucune racine rationnelle (cas des polynômes de Legendre).

Il n'existe pas à notre connaissance de résultats généraux sur l'irrationalité des racines de polynômes de Jacobi. Dans \cite{Ren13}, Render établit des résultats partiels donnant seulement des indications sur la forme des racines rationnelles possibles dans les cas particuliers où $2\alpha$ et $2\beta$ sont entiers.

\section{Polynômes de Laguerre et d'Hermite}
\subsection{Polynômes de Laguerre}
\subsubsection{Définition}

\begin{definition} \hspace{1cm}

\begin{enumerate}
\item Pour tout $d\in\N$, on définit le polynôme de Laguerre d'indice $d$ par
\[L_d:= \sum_{k=0}^{d} \frac{(-1)^k}{k!}\binom{d}{k} X^k.\]
\item Si $\alpha>-1$ est un réel quelconque, on définit, pour tout $d\in\N$, le polynôme de Laguerre généralisé de paramètre $\alpha$ et d'indice $d$ par
\[L_{d,\alpha}:=\sum_{k=0}^{d} \dfrac{(-1)^k}{k!}\dbinom{\alpha+d}{d-k} X^k.\]
En particulier, pour tout $d\in\N$, $L_d=L_{d,0}$.
\end{enumerate} 
\end{definition}

\begin{rem} --- On trouve souvent la notation $L_d^{(\alpha)}$ pour désigner le polynôme de Laguerre généralisé de paramètre $\alpha$ et d'indice $d$ mais nous avons préféré utiliser la notation $L_{d,\alpha}$ pour éviter la confusion, lorsque $\alpha\in\N$, avec la dérivée d'ordre $\alpha$ de $L_d$.
\label{rem_notation_Laguerre}
\end{rem}

\begin{propriete} --- Soit $\alpha>-1$ et $w_{\alpha}: t \mapsto t^{\alpha}\mathrm{e}^{-t}$ définie sur $I=\intervalleoo{0}{+\infty}$. La famille $(L_{n,\alpha})$ est une famille de polynômes $w_{\alpha}-orthogonaux$.
\label{prop_Laguerre_orthogonaux}
\end{propriete}

\emph{Preuve}. --- Pour tout polynôme $P\in \R[X]$, $P(t)w_{\alpha}(t)=\mathop{O}\limits_{t\to 0}(t^{\alpha})$ donc la fonction $t\mapsto P(t)w_{\alpha}(t)$ est intégrable au voisinage de $0$ car $\alpha>-1$ et $P(t)w_{\alpha}(t) =\mathop{o}\limits_{t\to+\infty}\left(\frac{1}{t^2}\right)$ donc la fonction $t\mapsto P(t)w_{\alpha}(t)$  au voisinage de $+\infty$. Comme, par ailleurs, $w_{\alpha}$ est clairement continue, positive et non identiquement nulle, on conclut que $w_{\alpha}$ est une fonction poids sur $\intervalleoo{0}{+\infty}$.

De plus, il est immédiat que, pour tout $n\in\N$, $L_{n,\alpha}$ est un polynôme de degré $n$.

Soit $d\in\N$. Considérons la fonction $f_{d,\alpha}:t \mapsto \frac{t^{\alpha+d}}{d!}\mathrm{e}^{-t}$ définie sur $I$. Elle est clairement de classe $C^{\infty}$ sur $I$ et, d'après la formule de Leibniz, pour tout $j\in\llbracket 0, d\rrbracket$ et tout $t\in I$,
\[f_{d,\alpha}^{(j)}(t)=\sum_{k=0}^j \binom{j}{k} (-1)^k\mathrm{e}^{-t}\dfrac{(\alpha+d-(j-k)+1)_{j-k}}{d!}t^{\alpha+d-(j-k)} = G_{d,j,\alpha}(t)t^{\alpha}\mathrm{e}^{-t} \]
en posant, pour tout $j\in\llbracket 0, d\rrbracket$, 
\[G_{d,j,\alpha}:=\sum_{k=0}^j \binom{j}{k} (-1)^k\dfrac{(\alpha+d+k-j+1)_{j-k}}{d!}X^{d+k-j}.\]
De plus, pour tout $j\in\llbracket 0, d\rrbracket$ et tout $k\in\llbracket 0, j\rrbracket$, $0 \leqslant d+k-j \leqslant d$ donc $G_{d,j,\alpha}\in \R_d[X]$.

Si $j\in\llbracket 0, d-1\rrbracket$, la valuation de $G_{d,j,\alpha}$ est $d-j>0$. De plus, par définition, 
\[G_{d,d,\alpha}=\sum_{k=0}^d \binom{d}{k} (-1)^k\dfrac{(\alpha+k+1)_{d-k}}{d!}X^{k}=\sum_{k=0}^d \dfrac{(-1)^k}{k!}\dbinom{\alpha+d}{d-k}X^{k}=L_{d,\alpha}.\]

Soit $m$ et $n$ deux entiers naturels tels que $m < n$. Alors, en intégrant par parties,
\begin{align*}
\varphi_{w_{\alpha}}(L_{m,\alpha}, L_{n,\alpha}) &= \int_{0}^{+\infty} L_{m,\alpha}(t)L_{n,\alpha}(t)t^{\alpha}\mathrm{e}^{-t} \diff t \\
& = \int_{0}^{+\infty} L_{m,\alpha}(t)f_{n,\alpha}^{(n)}(t) \diff t \\
&= \left[L_{m,\alpha}(t)f_{n,\alpha}^{(n-1)}(t)\right]_0^{+\infty} -  \int_{0}^{+\infty} L_{m,\alpha}'(t)f_{n,\alpha}^{(n-1)}(t) \diff t \\
&=\left[L_{m,\alpha}(t)G_{n,n-1,\alpha}(t)t^{\alpha}\mathrm{e}^{-t}\right]_0^{+\infty} -  \int_{0}^{+\infty} L_{m,\alpha}'(t)f_{n,\alpha}^{(n-1)}(t) \diff t \\
&=-  \int_{0}^{+\infty} L_{m,\alpha}'(t)f_{n,\alpha}^{(n-1)}(t) \diff t
\end{align*}
car la valuation de $G_{n,n-1,\alpha}$ est $1$ et $\alpha>-1$ donc $G_{n,n-1,\alpha}(t)t^{\alpha} \xrightarrow[t\to 0]{} 0$. En itérant les intégrations par parties, on aboutit à
\begin{align*}
\varphi_{w_{\alpha}}(L_{m,\alpha}, L_{n,\alpha})&=(-1)^m \int_{0}^{+\infty} L_{m,\alpha}^{(m)}(t)f_{n,\alpha}^{(n-m)}(t) \diff t = \int_{0}^{+\infty} f_{n,\alpha}^{(n-m)}(t) \diff t \\
&=\left[f_{n,\alpha}^{(n-m-1)}(t)\right]_0^{+\infty} = \left[G_{n,n-m-1,\alpha}(t)t^{\alpha}\mathrm{e}^{-t}\right]_0^{+\infty}=0
\end{align*}
car la valuation de $G_{n,n-m-1,\alpha}$ est $m+1\geqslant 1$.

Ainsi, $(L_{n,\alpha})$ est bien une famille de polynômes $w_{\alpha}-$orthogonaux. \hfill $\square$

\bigskip

On déduit alors de la propriété \ref{prop_poly_ortho_racines_simples} le résultat suivant.

\begin{corollaire} --- Soit un réel $\alpha>-1$. Pour tout entier $d\geqslant 1$, $L_{d,\alpha}$ est scindé à racines simples sur $\R$ et toutes ses racines appartiennent à $\intervalleoo{0}{+\infty}$.
\label{coro_Laguerre_racines_simples}
\end{corollaire}

\subsubsection{Irrationalité des racines}

\begin{theoreme} --- Pour tout entier $d\geqslant 2$, les $d$ racines réelles de $L_d$ sont irrationnelles.
\label{theo_Laguerre_0}
\end{theoreme}

\emph{Preuve}. --- Soit un entier $d\geqslant 2$. Posons, pour tout $k\in\llbracket 0, d \rrbracket$, $a_k=(-1)^k\binom{d}{k} \in \Z$ et $b_k=\frac{d!}{k!} \in \Z$. Alors, par définition, $d!L_d=\sum\limits_{k=0}^d a_kb_k X^k$. Soit $p$ un diviseur premier de $d$. Pour tout $k\in\llbracket 0, d-1\rrbracket$, $b_k=d(d-1)\cdots(k+1)$ donc $p$ divise $b_k$. De plus, d'après le corollaire \ref{coro_Legendre}, pour tout $k\in\llbracket 1, d\rrbracket$,
\[v_p(b_k)=v_p(d!)-v_p(k!) > v_p(b_0)-\dfrac{k}{p-1} > v_p(b_0)-k.\]
Par ailleurs, $\abs{a_0}=\abs{a_d}=b_d=1$ donc $p$ ne divise ni $a_0$, ni $a_d$, ni $b_d$. On déduit donc de la propriété \ref{prop_cri_irra_poly_Newton} que les racines réelles de $d!L_d$ sont irrationnelles et donc celle de $L_d$ aussi. \hfill $\square$

\bigskip

Si $\alpha>-1$ est un rationnel quelconque, le polynôme $L_{d,\alpha}$ peut avoir des racines rationnelles. Ainsi, on vérifie que $6$ est racine de $L_{2,2}$, $L_{2,7}$ et $L_{4,5}$, que $\frac{15}{2}$ est racine de $L_{3,\frac{21}{2}}$ ou encore que $\frac{12}{5}$ est racine de $L_{5,\frac{7}{5}}$. Cependant, ceci ne peut se produire si $d$ est \og suffisamment grand \fg{}.

\begin{theoreme}[Filaseta \& Lam, 2002] --- Soit un rationnel $\alpha>-1$. Il existe un entier $M_{\alpha}$ tel que, pour tout $d\geqslant M_{\alpha}$, les $d$ racines réelles de $L_{d,\alpha}$ sont toutes irrationnelles.
\label{theo_Laguerre_généralisé}
\end{theoreme}

La démonstration de ce théorème va nécessiter le lemme suivant dû à A. Thue \cite{Thu08} (voir également \cite{Pol18}).

\begin{lemme} --- Soit $a$, $b$, $c$ et $d$ des entiers relatifs tels que $ac\neq 0$ et $bc-ad \neq 0$. On note $n_0$ le plus petit entier naturel tel que, pour tout $n\geqslant n_0$, $\abs{(an+b)(cn+d)}\geqslant 2$ et, pour tout $n\geqslant n_0$, on désigne par $u_n$ le plus grand diviseur premier de $(an+b)(cn+d)$. Alors, la suite $(u_n)_{n\geqslant n_0}$ diverge vers $+\infty$.
\label{lem_Thue}
\end{lemme}

Nous démontrons ce lemme en admettant le théorème suivant lui aussi prouvé, dans un cas plus général, par A. Thue \cite{Thu09} (voir aussi \cite{AJ89} et \cite[p. 177 et suiv.]{ST15} pour une présentation exhaustive du cas particulier suivant).

\begin{theoreme}[Thue, 1909] --- Soit $u$ et $v$ deux entiers relatifs non nuls. Pour tout entier relatif non nul $w$, l'équation $ux^3-vy^3=w$ admet un nombre fini de solutions $\couple{x}{y}\in\Z^2$.
\label{theo_Thue}
\end{theoreme}

\noindent \textbf{\emph{Démonstration du lemme \ref{lem_Thue}}} --- Remarquons que, comme $ac\neq 0$, $\lim\limits_{n\to+\infty} \abs{(an+b)(cn+d)}=+\infty$ donc $n_0$ existe et $(u_n)_{n\geqslant n_0}$ est bien définie. Supposons, par l'absurde que $(u_n)$ ne tende pas vers $+\infty$. Alors, il existe une constante $K\geqslant 2$ telle que l'ensemble $E:=\{n\geqslant n_0 \mid u_n \leqslant K\}$ est infini. Notons $p_1$, $p_2$, ..., $p_m$ les nombres premiers inférieurs ou égaux à $K$. Alors, pour tout $n\in E$, les décompositions en produit de facteurs premiers de $an+b$ et de $cn+d$ ne contiennent que les nombres $p_1$, $p_2$, ..., $p_m$. Ainsi, pour tout $n\in E$, il existe des entiers naturels $\alpha_{1,n}$, ..., $\alpha_{m,n}$, $\beta_{1,n}$, ..., $\alpha_{m,n}$ tels que $an+b=\prod\limits_{i=1}^{m} p_i^{\alpha_{i,n}}$ et $cn+d=\prod\limits_{i=1}^{m} p_i^{\beta_{i,n}}$. Notons, pour tout $n\in E$ et tout $i\in\llbracket 1, m \rrbracket$, $r_{i,n}$ le reste de $\alpha_{i,n}$ modulo 3 et $s_{i,n}$ celui de $\beta_{i,n}$. Il existe alors des entiers $x_n$ et $y_n$ tels que $an+b=\left(\prod\limits_{i=1}^m p_i^{r_{i,n}}\right)x_n^3$ et $cn+d=\left(\prod\limits_{i=1}^m p_i^{s_{i,n}}\right)y_n^3$. Pour tout $n\in E$ et tout $i\in\llbracket 1, m\rrbracket$, $r_{i,n}\in\{0, 1, 2\}$ et $s_{i,n}\in\{0, 1, 2\}$ donc le $2m-$uplet $(r_{1,n}, ..., r_{m,n}, s_{1,n}, ..., s_{m,n})$ ne peut prendre que $3^{2m}$ valeurs distinctes. Comme $E$ est infini, il existe donc un $2m-$uplet $(r_1, ..., r_m, s_1, ..., s_m)$ et un ensemble infini $F\subset E$ tels que, pour tout $n\in F$ et tout $i\in\llbracket 1, m\rrbracket$, $r_{i,n}=r_i$ et $s_{i,m}=s_i$. Posons alors $g:=\prod\limits_{i=1}^m p_i^{r_i}$, $u:=cg \neq 0$, $h:=\prod\limits_{i=1}^m p_i^{s_i}$, $v:=ah\neq 0$ et $w:=bc-ad \neq 0$. Pour tout $n\in F$,
\[w=c(an+b)-a(cn+d)=c\left(\prod_{i=1}^m p_i^{r_{i}}\right)x_n^3-a\left(\prod_{i=1}^m p_i^{s_{i}}\right)y_n^3=ux_n^3-vy_n^3.\]
Or, pour tout $\couple{n}{n'}\in F^2$ tel que $n\neq n'$, $an+b \neq a'n+b$ car $a\neq 0$ donc $x_n^3 = \frac{an+b}{g} \neq \frac{an'+b}{g}=x_{n'}^3$ et ainsi $x_n \neq x_{n'}$. Puisque $F$ est infini, on en déduit que l'équation $ux^3-vy^3=w$ admet une infinité de solutions dans  $\Z^2$, ce qui entre en contradiction avec le théorème \ref{theo_Thue}. Ainsi, $(u_n)_{n\geqslant n_0}$ tend vers $+\infty$ lorsque $n$ tend vers $+\infty$. \hfill $\square$

\bigskip

\emph{Preuve du théorème \ref{theo_Laguerre_généralisé}}. --- Si $\alpha=0$, le résultat découle du théorème \ref{theo_Laguerre_0} avec $M_0=2$. Supposons, dans la suite, $\alpha \neq 0$ et écrivons $\alpha$ sous forme irréductible $\alpha=\frac{s}{t}$ avec $s\neq 0$. Posons, pour tout $d\in\N$, $S_{d,\alpha}:=d! t^d L_{d,\alpha}$. Alors, pour tout entier $d\geqslant 2$,
\begin{align*}
S_{d,\alpha}&=\sum_{k=0}^d \dfrac{(-1)^k d!t^d\left(\frac{s}{t}+d\right)\left(\frac{s}{t}+d-1\right)\cdots\left(\frac{s}{t}+k+1\right)}{k!(d-k)!}X^k \\
&=\sum_{k=0}^d (-1)^k \binom{d}{k} t^{k} (s+dt)(s+(d-1)t)\cdots(s+(k+1)t)X^k.
\end{align*}
Posons, pour tout $k\in\llbracket 0, d \rrbracket$, $a_k=(-1)^k t^{k} \in\Z$ et $b_k=\binom{d}{k}(s+dt)(s+(d-1)t)\cdots(s+(k+1)t)\in \Z$ de telle sorte que $S_{d,\alpha}=\sum\limits_{k=0}^d a_kb_kX^k$.

Comme $s\neq 0$, d'après le lemme \ref{lem_Thue}, il existe un entier $M_{\alpha}$ tel que, pour tout $d\geqslant M_{\alpha}$, $d(s+dt)$ admet un diviseur premier $p>\abs{s}+t$. En particulier, $p>t$ donc $p$ ne divise pas $t$. Par suite, $p$ ne divise ni $a_0=1$ ni $a_d=(-1)^dt^d$.

\emph{1\up{er} cas}. -- Supposons que $p$ divise $s+dt$. Alors, $p$ divise $b_k$ pour tout $k\in\llbracket 0, d-1\rrbracket$ et $p$ ne divise pas $b_d=1$. De plus, pour tout $k\in\llbracket 1, d\rrbracket$,
\begin{align*}
v_p(b_k)&=v_p\left(\dbinom{d}{k}\dfrac{(s+dt)(s+(d-1)t)\cdots(s+t)}{(s+kt)(s+(k-1)t)\cdots(s+t)}\right)&&\\
&\geqslant v_p\left(\dfrac{(s+dt)(s+(d-1)t)\cdots(s+t)}{(s+kt)(s+(k-1)t)\cdots(s+t)}\right) &&\\
&\geqslant v_p(b_0)-v_p((s+kt)(s+(k-1)t)\cdots(s+t)) &&\\
& \geqslant v_p(b_0)-v_p((\abs{s}+kt)!)& & (\text{d'après le lemme \ref{lem_majoration_factorielle}})&\\
& > v_p(b_0) - \dfrac{\abs{s}+kt}{p-1} & & (\text{d'après le corollaire \ref{coro_Legendre}}) &\\
& > v_p(b_0) - \dfrac{\abs{s}+kt}{\abs{s}+t} & & (\text{car } p>\abs{s}+t)&
\end{align*}
Or, si $k\geqslant 1$, $k(\abs{s}+t)\geqslant \abs{s}+kt$ donc $\frac{\abs{s}+kt}{\abs{s}+t}\leqslant k$. Ainsi, pour tout $k\in\llbracket 1, d\rrbracket$, $v_p(b_k)>v_p(b_0)-k$ et on conclut, grâce à la propriété \ref{prop_cri_irra_poly_Newton}, que $S_{d,\alpha}$ n'a pas de racine rationnelle et donc $L_{d,\alpha}$ non plus.

\emph{2\up{ième} cas}. -- Supposons que $p$ divise $d$. Soit $k\in\llbracket 0, d-1\rrbracket$. Si $k>d-p$ alors, étant donné que $(d-k)!\binom{d}{k}=d\frac{(d-1)!}{k!}$, $p$ divise $(d-k)!\binom{d}{k}$ mais, comme $p>d-k$, $p$ est premier avec $(d-k)!$ donc, par le lemme de Gauss, $p$ divise $\binom{d}{k}$ et, par suite, $p$ divise $b_k$. Sinon, $k\leqslant d-p$ donc $d-k\geqslant p$. Notons, pour tout $u\in\N$, $\overline{u}$ la classe de $u$ dans $\EQ{\Z}{p\Z}$. Comme $\llbracket k+1, d\rrbracket$ contient au moins $p$ entiers consécutifs, $\{\overline{u} \mid u\in\llbracket k+1, d\rrbracket\}=\EQ{\Z}{p\Z}$. Or, comme $t$ est premier avec $p$, l'application $\overline{x}\mapsto \overline{s+tx}$ est une bijection de $\EQ{\Z}{p\Z}$ dans lui-même donc il existe $u\in\llbracket k+1, d\rrbracket$ tel que $\overline{s+tu}=\overline{0}$. Autrement dit, il existe un entier $u\in\llbracket k+1, d \rrbracket$ tel que $p$ divise $s+tu$ et, par suite, $p$ divise $b_k$. Ainsi, dans tous les cas, $p$ divise $b_k$. De plus, comme précédemment, pour tout $k\in\llbracket 1, d\rrbracket$, $v_p(b_k)>v_p(b_0)-k$ et on conclut de même. \hfill $\square$

\begin{lemme} --- Soit $\varepsilon\in\{-1, 1\}$. L'unique solution de l'équation $2^{m+1}+\varepsilon=3^q$ d'inconnue $\couple{m}{q}\in(\N^*)^2$ est $\couple{1}{1}$ si $\varepsilon=-1$ et $\couple{2}{2}$ si $\varepsilon=1$.
\label{lem_puissances_2_et_3}
\end{lemme}

\emph{Preuve}. --- Supposons que $\couple{m}{q}\in(\N^*)^2$ vérifie $2^{m+1}+\varepsilon=3^q$. 

Si $m=1$ alors $3^q=4+\varepsilon$ donc $\varepsilon=-1$ et $q=1$.

Supposons $m\geqslant 2$. Alors, $2^{m+1} \equiv 0 \pmod{8}$ donc $2^{m+1}-1 \equiv 7 \pmod{8}$ et $2^{m+1}+1 \equiv 1 \pmod{8}$. Or, comme $3^2\equiv 1 \pmod{8}$, $3^q\equiv 1 \pmod{8}$ si $q$ est pair et $3^q\equiv 3 \pmod{8}$ si $q$ est impair donc $\varepsilon=1$ et $q$ est pair. Il existe donc un entier $g\geqslant 1$ tel que $2^{m+1}+1=3^{2g}$ et ainsi $(3^g-1)(3^g+1)=2^{m+1}$. Par suite, il existe deux entiers $i$ et $j$ avec $1 \leqslant i < j$ tels que $3^g-1=2^i$ et $3^g+1=2^j$ donc $2^i=2^j-2=2(2^{j-1}-1)$ ce qui impose $i=1$ et $j=2$. On en déduit que $g=1$ donc $q=2$ et $m=2$. 

La réciproque est immédiate. \hfill $\square$

\begin{corollaire} --- Soit $\alpha\in\left\{-\frac{1}{2}, \frac{1}{2}, 1\right\}$. Pour tout entier $d\geqslant 2$, les $d$ racines réelles de $L_{d,\alpha}$ sont irrationnelles.
\label{coro_Laguerre}
\end{corollaire}

\emph{Preuve}. --- La démonstration du théorème \ref{theo_Laguerre_généralisé} montre que, si l'écriture sous forme de fraction irréductible de $\alpha$ est $\frac{s}{t}$ alors $M_{\alpha}$ est n'importe quel entier tel qu'il existe un nombre premier $p\geqslant M_{\alpha}$ tel que $p$ divise $d(s+dt)$ et $p>\abs{s}+t$. Il suffit donc de montrer que si $d\geqslant 2$ alors $d(s+dt)$ admet un diviseur premier $p>\abs{s}+t$.

Soit un entier $d\geqslant 2$.

Si $\alpha=1$ alors $s=t=1$. Dès lors, $d(s+dt)=d(d+1)$ est le produit de deux entiers consécutifs donc l'un d'eux est un nombre impair au moins égal à 3 car $d\geqslant 2$. Il s'ensuit que $d(d+1)$ admet un diviseur premier impair $p \geqslant 3 > \abs{s}+t=2$. 

Si $\abs{\alpha}=\frac{1}{2}$ alors $\abs{s}=1$ et $t=2$ et donc $d(s+dt)=d(2d+s)$ avec $s\in\{-1,1\}$. Comme $\abs{s}+t=3$, si $d(2d+s)$ admet un diviseur premier $p\geqslant 5$, le résultat est acquis. Dans le cas contraire, les seuls diviseurs premiers de $d$ et $2d+s$ sont $2$ et $3$. Or, $d \geqslant 2$ et $2d+s\geqslant 3$ sont premiers entre eux (car $\abs{s}=1$) et $2d+s$ est impair donc il existe deux entiers $m \geqslant 1$ et $q\geqslant 1$ tels que $d=2^{m}$ et $2d+s=3^q$ et, ainsi, $2^{m+1}+s=3^q$. On déduit alors du lemme \ref{lem_puissances_2_et_3} que $s=-1$ et $m=1$ ou $s=1$ et $m=2$ i.e. $\alpha=-\frac{1}{2}$ et $d=2$ ou $\alpha=\frac{1}{2}$ et $d=4$. Il suffit donc, pour conclure, de montrer que les racines de $L_{2,-\frac{1}{2}}$ et $L_{4,\frac{1}{2}}$ sont irrationnelles. 

D'une part, les racines de $L_{2,-\frac{1}{2}}=\frac{1}{2}X^2-\frac{3}{2}X+\frac{3}{8}$ sont $\frac{3-\sqrt{6}}{2}$ et $\frac{3+\sqrt{6}}{2}$ donc elles sont bien irrationnelles. 

D'autre part, posons $P:=2^44!L_{4,\frac{1}{2}}=16X^4-288X^3+1512X^2-2520X+945$. On peut utiliser la propriété \ref{prop_RRT} pour montrer que les racines de $P$ (et donc de $L_{4,\frac{1}{2}}$) sont irrationnelles mais c'est assez fastidieux car $16$ et $945$ ont beaucoup de diviseurs (il y a en tout 80 valeurs à tester). Nous pouvons améliorer les choses en utilisant une autre méthode due à Newton (voir \cite{MS01}). Supposons que $P$ admette une racine rationnelle écrite sous forme irréductible $\frac{u}{v}$. Alors, $vX-u$ divise $P$ dans $\Q[X]$ et donc, d'après le lemme de Gauss (\cite[§ 42, p. 26-27]{Gau01}), $vX-u$ divise $P$ dans $\Z[X]$. Ainsi, on peut écrire $P=(vX-u)Q$ avec $Q\in\Z[X]$. On a alors $P(-2)=(-2v-u)Q(-2)$ et $P(-1)=(-v-u)Q(-1)$ donc les entiers $P(-2)$ et  $P(-1)$ admettent des diviseurs respectifs $\delta_{-2}=-2v-u$ et $\delta_{-1}=-v-u$ qui diffèrent de $v$. Or, $P(-2)=14593$ est premier donc $\delta_{-2}\in\{-14593, -1, 1, 14593\}$. De plus, $v>0$ et, d'après la proposition \ref{prop_RRT}, $v$ divise $16$ donc $v\in\{1, 2, 4, 8, 16\}$. Or, $P(-1)=5281$ est également premier donc la seule valeur possible pour $\delta_{-1}$ est $1$ obtenue avec $v=2$. Par suite, $u=-v-\delta_{-1}=-3$ ce qui est absurde car $\frac{u}{v}>0$ d'après le corollaire \ref{coro_Laguerre_racines_simples} donc les racines de $P$ sont irrationnelles. \hfill $\square$

\bigskip

\begin{rem} --- L'irréductibilité des polynômes de Laguerre $L_d$ pour tout $d\geqslant 2$ a été prouvée par Schur en 1929 \cite{Sch29} en s'appuyant sur une généralisation du postulat de Bertrand qui avait été déjà été obtenue par Sylvester \cite{Syl92} et que Schur a redémontré indépendamment. Dans l'article \cite{FL02}, Filaseta et Lam ont en fait montré que, pour tout rationnel $\alpha$, il existe un entier $M_{\alpha}$ tel que, pour tout $d\geqslant M_{\alpha}$, le polynôme $L_{d,\alpha}$ est irréductible sur $\Q$. Leur démonstration utilise le lemme \ref{lem_Thue}, la propriété énoncée dans la remarque \ref{rem_poly_Newton} ainsi que des estimations sur la répartition des nombres premiers dans une progression arithmétique donnée.
\label{rem_Laguerre_irré}
\end{rem}

\subsection{Polynômes d'Hermite}

\begin{definition} --- Pour tout $d\in\N$, on définit le polynôme d'Hermite d'indice $d$ par
\[H_d:= \sum_{k=0}^{\entiere{\frac{d}{2}}} \left(-\frac{1}{2}\right)^k\dfrac{d!}{(d-2k)!k!}X^{d-2k}.\]
\end{definition} 

On montre sans difficulté que, pour tout $d\in\N$, $H_d$ est un polynôme unitaire de degré $d$ de $\Q[X]$. En remarquant que, pour tout $d\in\N$,
\[H_d:= \sum_{k=0}^{\entiere{\frac{d}{2}}} (-1)^k\binom{d}{2k}\prod_{j=1}^k (2j-1) X^{d-2k}\]
on peut affirmer que $H_d\in\Z[X]$. De plus, si on note $w : t\mapsto \mathrm{e}^{-\frac{t^2}{2}}$ définie sur $\R$, alors on vérifie en raisonnant par récurrence et selon la parité de $d$ que, pour tout réel $t$, $w^{(d)}(t)=(-1)^dH_d(t)f_d(t)$. De là, on déduit, comme pour la propriété \ref{prop_Laguerre_orthogonaux}, que $(H_n)$ est une famille de polynômes $w-$orthogonaux et ainsi la propriété \ref{prop_poly_ortho_racines_simples} assure que, pour tout $d\in\N^*$, le polynôme $H_d$ est scindé à racines simples sur $\R$. 

\bigskip

On peut cependant ramener l'étude des polynômes d'Hermite à celle des polynômes de Laguerre généralisés grâce au lemme suivant.

\begin{lemme} --- Soit $n\in\N$. Alors, 
\[H_{2n}=(-2)^{n} n! L_{n,-\frac{1}{2}}\left(\dfrac{X^2}{2}\right) \qquad \text{et} \qquad H_{2n+1}=(-2)^{n} n! XL_{n,\frac{1}{2}}\left(\frac{X^2}{2}\right).\]
\label{lem_Hermite_Laguerre}
\end{lemme}

\emph{Preuve}. --- En utilisant le lemme \ref{lem_coeff_binomial},
\begin{align*}
H_{2n}&=\sum_{k=0}^n \dfrac{(-1)^k}{2^k}\cdot\dfrac{(2n)!}{(2n-2k)!k!}X^{2n-2k} =\sum_{k=0}^n \dfrac{(-1)^{n-k}}{2^{n-k}}\cdot\dfrac{(2n)!}{(2k)!(n-k)!}\left(X^2\right)^k \\
&=(-1)^n\sum_{k=0}^n \dfrac{(-1)^k}{k!} 2^{n-k}n!\dfrac{k!(2n)!}{4^{n-k}(n-k)!n!(2k)!}\left(X^2\right)^k \\
&=(-2)^nn! \sum_{k=0}^n \dfrac{(-1)^k}{k!} \dbinom{-\frac{1}{2}+n}{n-k}\left(\dfrac{X^2}{2}\right)^k =(-2)^nn!L_{n,-\frac{1}{2}}\left(\dfrac{X^2}{2}\right)
\end{align*}
et
\begin{align*}
H_{2n+1}&=\sum_{k=0}^n \dfrac{(-1)^k}{2^k}\cdot\dfrac{(2n+1)!}{(2n+1-2k)!k!}X^{2n+1-2k} =X\sum_{k=0}^n \dfrac{(-1)^{n-k}}{2^{n-k}}\cdot\dfrac{(2n+1)!}{(2k+1)!(n-k)!}\left(X^2\right)^k \\
&=(-1)^nX\sum_{k=0}^n \dfrac{(-1)^k}{k!} 2^{n-k}n!\dfrac{2(k+1)(n+1)k!(2n+1)!}{2(k+1)(n+1)4^{n-k}(n-k)!(2k+1)!n!}\left(X^2\right)^k \\
&=(-1)^nX\sum_{k=0}^n \dfrac{(-1)^k}{k!} 2^{n-k}n!\dfrac{(k+1)!(2n+2)!}{4^{n-k}(n-k)!(n+1)!(2k+2)!}\left(X^2\right)^k \\
&=(-2)^nn! X\sum_{k=0}^n \dfrac{(-1)^k}{k!} \dbinom{\frac{1}{2}+n}{n-k}\left(\dfrac{X^2}{2}\right)^k 
=(-2)^nn!XL_{n,\frac{1}{2}}\left(\dfrac{X^2}{2}\right).
\end{align*}
\hfill $\square$

\begin{theoreme} --- Pour tout entier $d\geqslant 1$, $H_d$ admet $d$ racines simples. De plus, si $d\geqslant 4$ est pair, les $d$ racines de $H_d$ sont toutes irrationnelles et, si $d\geqslant 3$ est impair, les $d-1$ racines non nulles de $H_d$ sont toutes irrationnelles.
\label{theo_Hermite}
\end{theoreme}

\emph{Preuve}. --- Soit un entier $d\geqslant 1$. 

Si $d$ est pair, il existe un entier $n\geqslant 1$ tel que $d=2n$. Alors, $H_{d}=H_{2n}=(-2)^nn!L_{n,-\frac{1}{2}}\left(\frac{X^2}{2}\right)$. D'après le corollaire \ref{coro_Laguerre_racines_simples}, $L_{n,-\frac{1}{2}}$ admet $n$ racines simples qui appartiennent toutes à $\intervalleoo{0}{+\infty}$. Or, si $r>0$ est une racine de $L_{n,-\frac{1}{2}}$ alors $\sqrt{2r}$ et $-\sqrt{2r}$ sont deux racines distinctes de $L_{n,-\frac{1}{2}}\left(\frac{X^2}{2}\right)$ donc $H_d$ admet $2n=d$ racines simples. De plus, si $r\in\Q$ alors $\frac{r^2}{2}\in\Q$ et, pour $n\geqslant 2$, $L_{n,-\frac{1}{2}}$ n'a pas de racine rationnelle d'après le corollaire \ref{coro_Laguerre}, donc, si $d\geqslant 4$, $H_d$ n'a pas de racines rationnelles. En revanche, $H_2=X^2-1=(X-1)(X+1)$ possède des racines rationnelles. 

Si $d$ est impair, on écrit $d=2n+1$ avec $n\in\N$. Alors, $H_{d}=H_{2n+1}=(-2)^nn!XL_{n,\frac{1}{2}}\left(\frac{X^2}{2}\right)$. D'après le corollaire \ref{coro_Laguerre_racines_simples}, $L_{n,\frac{1}{2}}$ admet $n$ racines simples qui appartiennent toutes à $\intervalleoo{0}{+\infty}$ donc, comme précédemment, $L_{n,\frac{1}{2}}\left(\frac{X^2}{2}\right)$ admet $2n$ racines simples strictement positives et ainsi $H_d$, qui s'annule également en $0$, possède $2n+1=d$ racines simples. De plus, comme précédemment, $L_{n,\frac{1}{2}}$ n'a pas de racines rationnelles si $n\geqslant 2$ donc $H_d$ n'a pas de racines rationnelles autres que $0$ si $d\geqslant 5$. Par ailleurs, $H_3=X^3-3X=X(X-\sqrt{3})(X+\sqrt{3})$ donc cette dernière affirmation reste vraie pour $d=3$. \hfill $\square$

\bigskip

\begin{rem} --- Il est possible de démontrer directement l'irrationalité des racines de $H_d$ pour tout entier $d\geqslant 3$, sans faire référence aux polynômes de Laguerre, en remarquant que, pour tout $n\geqslant 1$,
\[H_{2n}=\sum_{j=0}^{n} (-1)^{n-j} \dbinom{n}{j} \dfrac{u_{2n}}{u_{2j}} X^{2j} \quad \text{ et } \quad H_{2n+1}=X\sum_{j=0}^{n} (-1)^{n-j} \dbinom{n}{j} \dfrac{u_{2n+2}}{u_{2j+2}} X^{2j}\]
où, pour tout entier $j \geqslant 1$, $u_{2j}=1\times 3 \times 5 \times \cdots \times (2j-1) = \frac{(2j)!}{2^j j!}$
et en utilisant la propriété \ref{prop_cri_irra_poly_Newton}.

Par ailleurs, dans la deuxième partie de son article \cite{Sch29}, Schur a démontré l'irréductibilité de $H_d$ pour tout entier pair $d\geqslant 4$ et de $\frac{H_d}{X}$ pour tout entier impair $d\geqslant 3$.
\label{rem_autre_méthode_Hermite}
\end{rem}

\section{Polynômes de Bessel}

\subsection{Définition et propriétés}

\begin{definition} ---Pour tout entier $d\in\N$, on définit le polynôme de Bessel d'indice $d$ par
\[\mathscr{B}_d:=\sum_{k=0}^{d} \dfrac{(d+k)!}{2^kk!(d-k)!}X^k.\]
\end{definition}

Si $d\in\N$ alors, pour tout $k\in\llbracket 0, d \rrbracket$, le coefficient d'indice $k$ de $\mathscr{B}_d$ est 
\[\dfrac{(d+k)!}{2^k k! (d-k)!} = \dfrac{(d+k)!}{(2k)! (d-k)!}\times \dfrac{(2k)!}{2^k k!} = \dbinom{d+k}{2k}\times 1\times 3 \times \cdots \times (2k-1)\]
donc $\mathscr{B}_d$ est un polynôme de $\Z[X]$ de degré $d$. 

\begin{propriete} --- Soit un entier $d\geqslant 1$. Alors, $\mathscr{B}_d$ ne possède pas de racine réelle si $d$ est pair et $\mathscr{B}_d$ possède une unique racine réelle simple $\xi_d$ si $d$ est impair et, de plus, $\xi_d<0$.
\end{propriete}

\emph{Preuve}. --- Nous suivons ici la méthode exposée dans  \cite{Bur51}. 

Posons $Q_d:=X^d\mathscr{B}_d\left(\frac{1}{X}\right)$ qui est également un polynôme de $\Z[X]$ de degré $d$ et considérons les fonctions $f : x \mapsto \mathrm{e}^{-x}Q_d(x)$ et $g : x \mapsto f(-x)$ définies sur $\R$.

Montrons tout d'abord que $f$ et $g$ sont deux solutions sur $\R$ de l'équation différentielle
\[(E)~xy''-2dy'-xy=0.\]
Notons, pour tout $k\in\llbracket 0, d\rrbracket$, $c_k=\frac{(2d-k)!}{2^{d-k}(d-k)!k!}$ de telle sorte que $Q_d=\sum\limits_{k=0}^d c_k X^k$. Pour tout réel $x$, $f'(x)=\mathrm{e}^{-x}\left[-Q_d(x)+Q_d'(x)\right]$ et $f''(x)=\mathrm{e}^{-x}\left[Q_d(x)-2Q_d'(x)+Q_d''(x)\right]$ donc
\[xf''(x)-2df'(x)-xf(x)=\mathrm{e}^{-x} A(x)\]
avec
\begin{align*}
A(x)&=2dQ_d(x)-2(d+x)Q_d'(x)+xQ_d''(x) =\sum_{k=0}^d \left[2dc_kx^k -2(d+x)kc_kx^{k-1} + k(k-1)c_k x^{k-1}\right]\\
&=\sum_{k=0}^d 2(d-k)c_kx^k+ \sum_{k=0}^d k(k-1-2d)c_k x^{k-1}  =\sum_{k=0}^{d-1} 2(d-k)c_kx^k+ \sum_{k=1}^d k(k-1-2d)c_k x^{k-1}  \\
&=\sum_{k=0}^{d-1} 2(d-k)c_kx^k+ \sum_{k=0}^{d-1} (k+1)(k-2d)c_{k+1} x^{k} =\sum_{k=0}^{d-1} \left[2(d-k)c_k  + (k+1)(k-2d)c_{k+1} \right] x^{k} 
\end{align*}
Or, pour tout $k\in\llbracket 0, d-1\rrbracket$,
\begin{align*}
2(d-k)c_k  + (k+1)(k-2d)c_{k+1}&=\dfrac{2(d-k)(2d-k)!}{2^{d-k} (d-k)! k!} + \dfrac{(k+1)(k-2d)(2d-k-1)!}{2^{d-k-1} (d-k-1)! (k+1) !} \\
&=\dfrac{(2d-k)!}{2^{d-k-1} (d-k-1)! k!} - \dfrac{(2d-k)(2d-k-1)!}{2^{d-k-1} (d-k-1)! k!}\\
&=0
\end{align*}
Ainsi, pour tout réel $x$, $xf''(x)-2df'(x)-xf(x)=0$ donc $f$ est bien solution de $(E)$ sur $\R$.

En notant que, pour tout réel $x$, $-\left[(-x)f''(-x)+2d(-f'(-x))-(-x)f(-x)\right]=0$, on en déduit que, pour tout réel $x$, $xg''(x)-2dg'(x)-xg(x)=0$ donc $g$ est également solution de $(E)$ sur $\R$.

Remarquons que $f$ et $g$ sont linéairement indépendantes sur $\R$. Supposons, en effet, qu'il existe deux constantes réelles $\lambda$ et $\mu$ telles que, pour tout réel $x$, $\lambda f(x)+\mu g(x)=0$. Alors, pour tout réel $x$, $\lambda \mathrm{e}^{-x}Q_d(x)+\mu \mathrm{e}^x Q_d(-x)=0$. En utilisant le fait qu'un polynôme est négligeable devant la fonction exponentielle aux voisinages de $+\infty$ et $-\infty$, on en déduit que $\lim\limits_{x\to+\infty} \mu \mathrm{e}^x Q_d(-x)=0$ et $\lim\limits_{x\to-\infty} \lambda \mathrm{e}^{-x} Q_d(x)=0$. Or, $\lim\limits_{x\to+\infty} \abs{\mathrm{e}^x Q_d(-x)}=\lim\limits_{x\to-\infty} \abs{\mathrm{e}^{-x} Q_d(x)}=+\infty$ donc $\lambda=\mu=0$ et ainsi $f$ et $g$ sont linéairement indépendantes sur $\R$.

Ainsi, $(g, f)$ est un système fondamental de solutions de $(E)$. On en déduit, dès lors, que le wronskien $W:x\mapsto \begin{vmatrix} g(x) & f(x) \\ g'(x) & f'(x) \end{vmatrix}$ ne s'annule pas sur $\R$ et qu'il est, de plus, solution de l'équation différentielle $y'=\frac{2d}{x}y$ sur chacun des deux intervalles $\R_-^*$ et $\R_+^*$. Il s'ensuit qu'il existe deux constantes non nulles $K_1$ et $K_2$ telles que, pour tout $x<0$, $W(x)=K_1\mathrm{e}^{2d\ln (-x)}=K_1 x^{2d}$ et, pour tout $x>0$, $W(x)=K_2\mathrm{e}^{2d\ln x}=K_2 x^{2d}$. Comme $f $ et $g$ sont indéfiniment dérivables, $W$ est de classe $\mathcal{C}^{\infty}$ sur $\R$ donc $\lim\limits_{x\to 0^-} W^{(2d)}(x) = \lim\limits_{x\to 0^+} W^{(2d)}(x)$ i.e. $(2d)!K_1=(2d)!K_2$ donc $K_1=K_2$. Ainsi, il existe une constante $C_d\neq 0$ telle que, pour tout réel $x\neq 0$, $W(x)=C_dx^{2d}$ et, par continuité, pour tout réel $x$, $W(x)=C_dx^{2d}$ i.e. $f'(x)g(x)-f(x)g'(x)=C_dx^{2d}$.

En remplaçant $f$ et $g$ par leurs expressions en fonction de $Q_d$, on en déduit que, pour tout réel $x$, 
\[\mathrm{e}^{-x}\left[-Q_d(x)+Q_d'(x)\right]\mathrm{e}^{x}Q_d(-x)-\mathrm{e}^{-x}Q_d(x)\mathrm{e}^{x}\left[Q_d(-x)-Q_d'(-x)\right]=C_dx^{2d}\]
i.e.
\begin{equation}
Q_d'(x)Q_d(-x)+Q_d(x)Q_d'(-x)-2Q_d(x)Q_d(-x)=C_dx^{2d}.
\label{eq_diff_Bessel}
\end{equation}

Les coefficients de $\mathscr{B}_d$ (et donc de $Q_d$) sont des entiers naturels non nuls donc si $Q_d$ admet des racines réelles, elles sont nécessairement strictement négatives. Supposons que $\beta$ est une racine de $Q_d$. En substituant dans $\eqref{eq_diff_Bessel}$, on obtient $Q_d'(\beta)Q_d(-\beta)=C_d\beta^{2d}$. Comme $C_d$ n'est pas nulle, on en déduit que $Q_d'(\beta)$ n'est pas nul. De plus, comme $-\beta>0$ et comme les coefficients de $Q_d$ sont strictement positifs, $Q_d(-\beta)>0$. Ainsi, $\beta$ est une racine simple de $Q_d$ et $Q_d'(-\beta)$ est du signe de $C_d$. Supposons que $Q_d$ admette au moins deux racines réelles. Considérons alors deux racines consécutives $\beta_1 < \beta_2$ de $Q_d$ i.e. deux racines de $Q_d$ telles que $Q_d$ n'a pas de racine dans $\intervalleoo{\beta_1}{\beta_2}$. D'après ce qui précède, $Q_d'(\beta_1)$ et $Q_d'(\beta_2)$ sont non nuls et de même signe (celui de $C_d$) donc, sur des voisinages de $\beta_1$ et de $\beta_2$, la fonction $h : x\mapsto Q_d(x)$ a le même sens de variation. Quitte à considérer $-h$, on peut supposer que $h$ est croissante sur un voisinage de $\beta_1$ et sur un voisinage de $\beta_2$. Il existe alors un réel $\gamma_1>\beta_1$ et un réel $\gamma_2<\beta_2$ tels que $h(\gamma_1)>0$ et $h(\gamma_2)<0$. Par continuité de $h$, on en déduit que $h$ s'annule entre $\beta_1$ et $\beta_2$ ce qui est contradictoire avec l'hypothèse. Ainsi, $Q_d$ admet au plus une racine réelle (comptée avec multiplicité). 

Comme un réel non nul $\alpha$ est racine de $\mathscr{B}_d$ si et seulement si $\beta:=\frac{1}{\alpha}$ est racine de $Q_d$, on conclut que $\mathscr{B}_d$ admet également au plus une racine réelle (comptée avec multiplicité). Dès lors, si $d$ est pair, $Q_d$ n'a pas de racine réelle et si $d$ est impair, $d$ admet exactement une racine réelle simple d'après le théorème des valeurs intermédiaires. \hfill $\square$

\bigskip

\begin{rem} --- En déterminant un équivalent en $+\infty$ des deux membres de l'égalité \eqref{eq_diff_Bessel}, on obtient facilement que $C_d=2(-1)^{d+1}$ mais ceci n'est pas nécessaire pour la démonstration précédente.
\label{rem_Bessel_Cd}
\end{rem}

\subsection{Irrationalité des racines}

Nous allons démontrer que, pour tout entier impair $d\geqslant 3$, l'unique racine réelle de $\mathscr{B}_d$ est irrationelle. La démarche suivie reprend les idées de Filaseta et Trifonov exposées dans \cite{Fil95} et \cite{FT02}. Dans ces articles, les auteurs s'intéressent en fait à l'irréductibilité de $\mathscr{B}_d$, qui est démontrée dans $\cite{FT02}$ pour tout entier $d\geqslant 2$. Ils utilisent, pour cela, des estimations de certaines fonctions liées aux nombres premiers dont les majorations \og à la Tchebichef \fg{} suivantes, établies par Rosser et Schoenfeld \cite{RS62},
\[\forall x\geqslant 113{,}6 \quad \pi(x) \leqslant \dfrac{5x}{4\ln{x}} \qquad \text{et} \qquad \forall x>1  \quad \pi(x) \leqslant 1{,}25506\dfrac{x}{\ln{x}}\]
où $\pi(x)$ désigne le nombre de nombres premiers inférieurs ou égaux à $x$.

Pour notre preuve d'irrationalité, nous allons utiliser une majoration du même type qui est moins fine mais beaucoup plus élémentaire à démontrer.

\subsubsection{Une majoration \og à la Tchebichef \fg{}}

\begin{lemme} --- Soit un réel $c>0$. La fonction $T_c : x\mapsto \frac{x^c}{\ln x}$ est décroissante sur $\intervalleof{1}{\exp(c^{-1})}$ et croissante sur $\intervallefo{\exp(c^{-1})}{+\infty}$.
\label{lem_étude_fonction_Bessel}
\end{lemme}

\emph{Preuve}. --- Pour tout $x>1$, $T_c'(x)=\frac{x^{c-1}(c\ln x - 1)}{(\ln x)^2}$ donc $T_c'(x)$ est du signe de $c\ln x - 1$ et le résultat s'ensuit immédiatement. \hfill $\square$

\begin{propriete} --- Pour tout réel $x>1$, le nombre $\pi(x)$ de nombres premiers inférieurs ou égaux à $x$ vérifie $\pi(x) \leq \frac{(1+\ln 4)x}{\ln x}$.
\label{prop_majoration_pi}
\end{propriete}

\emph{Preuve}. --- Pour tout $n\in\N^*$, on note $\mathbb{P}_n$ l'ensemble des nombres premiers inférieurs ou égaux $n$ et $P_n=\prod_{p\in\mathbb{P}_n} p$. Remarquons que si $(p_k)_{k\in\N^*}$ désigne la suite des nombres premiers alors $P_n=\prod\limits_{k=1}^{\pi(n)} p_k$.

Soit $n\in\N^*$. \'Etant donné que $n!\binom{2n+1}{n}=\prod\limits_{j=2}^{n+1} (n+j)$, tout nombre premier $p\in\llbracket n+2, 2n+1 \rrbracket$ divise $n!\binom{2n+1}{n}$. Or, un tel nombre est premier avec $n!$ donc tout nombre premier $p\in\llbracket n+2, 2n+1 \rrbracket$ divise $\binom{2n+1}{n}$ et, par suite, $\binom{2n+1}{n}$ est divisible par $\frac{P_{2n+1}}{P_{n+1}}$. En particulier, $P_{2n+1} \leqslant \binom{2n+1}{n} P_{n+1}$. 

Par ailleurs, grâce à la formule du binôme de Newton,
\[\dbinom{2n+1}{n} = \dfrac{1}{2}\left[\dbinom{2n+1}{n}+\dbinom{2n+1}{n+1}\right] \leqslant \dfrac{1}{2}\sum_{j=0}^{2n+1} \dbinom{2n+1}{j} = 4^{n}\]
donc $P_{2n+1} \leqslant 4^{n}P_{n+1}$.

Montrons par récurrence forte que, pour tout $n\in\N^*$, $P_{n} \leqslant 4^{n-1}$. Comme $P_1=1$ et $P_2=2$, l'inégalité est vraie aux rangs $n=1$ et $n=2$. Supposons que, pour un certain entier $n\geqslant 2$, on ait $P_{j}\leqslant 4^{j-1}$ pour tout $j\in\llbracket 1, n\rrbracket$. Si $n$ est impair alors, comme $n+1$ est un entier pair différent de $2$, $P_{n+1}=P_n \leqslant 4^{n-1} \leqslant 4^{n}$. Si $n$ est pair alors il existe $k\in\N^*$ tel que $n=2k$ donc, d'après ce qui précède, 
\[P_{n+1}=P_{2k+1} \leqslant 4^{k}P_{k+1} \leqslant 4^{k}4^{k}=4^{2k} =4^{n}\]
donc l'inégalité est établie au rang $n+1$ dans tous les cas, ce qui permet de conclure.

Considérons, à présent, la suite $(u_n)$ définie pour tout $n\in\N^*$ par $u_{n}=\ln(n!)-n(\ln n-1)$. Alors, pour tout $n\in\N^*$, en utilisant le fait que, pour tout réel $x>-1$, $\ln(1+x)\leqslant x$,
\[u_{n+1}-u_n=\ln(n+1)-(n+1)(\ln(n+1)-1)+n(\ln n -1) =-n\ln\left(1+\dfrac{1}{n}\right)+1 \geqslant 0 \]
donc $(u_n)$ est croissante. En particulier, pour tout $n\in\N^*$, $u_n\geqslant u_1=1$ donc, pour tout $n\in\N^*$, $\ln(n!) \geqslant n(\ln n - 1)$.

Il s'ensuit que, pour tout entier $n\geqslant 2$, $\pi(n)(\ln(\pi(n))-1) \leqslant \ln(\pi(n)!)$. Or, pour tout $k\in\N^*$, $p_k\geqslant k$ donc $P_n\geqslant \prod\limits_{k=1}^{\pi(n)} k = \pi(n)!$ et ainsi $\pi(n)(\ln(\pi(n))-1) \leqslant \ln(P_n) \leqslant \ln (4^{n-1})$. On en déduit en particulier que $\pi(n)(\ln(\pi(n))-1) \leqslant (n-1)\ln 4 \leqslant n\ln 4$.

Considérons la fonction $f : x \mapsto \frac{1+\ln 4}{\ln x}\left[\ln\left(\frac{(1+ \ln 4)x}{\ln x}\right)-1\right]$ définie sur $\intervalleoo{1}{+\infty}$. Alors, pour tout réel $x>1$, $f'(x)=\frac{1+\ln 4}{x(\ln x)^2}\ln\left(\frac{\ln x}{1+\ln 4}\right)$
donc $f$ atteint son minimum en $x=4\mathrm{e}$ et ce minimum vaut $f(4\mathrm{e})=\ln 4$. On en déduit que, pour entier $n\geqslant 2$, 
\[\pi(n)(\ln(\pi(n))-1) \leqslant \dfrac{(1+\ln 4)n}{\ln n}\left[\ln\left(\dfrac{(1+\ln 4)n}{\ln n}\right)-1\right].\]
Remarquons que $x \mapsto x(\ln x - 1)$ est une primitive de $\ln$ sur $\intervalleoo{0}{+\infty}$ donc elle est strictement croissante sur $\intervallefo{1}{+\infty}$ et, ainsi, pour tout entier $n\geqslant 2$, $\pi(n) \leqslant \frac{(1+\ln 4)n}{\ln n}$.

D'après le lemme \ref{lem_étude_fonction_Bessel}, la fonction $g : x\mapsto \frac{(1+\ln 4)x}{\ln x}$ est décroissante sur $\intervalleof{1}{\mathrm{e}}$ et croissante sur $\intervallefo{\mathrm{e}}{+\infty}$: elle atteint son minimum en $x=\mathrm{e}$ et ce minimum est $g(\mathrm{e})=(1+\ln 4)\mathrm{e} > 6{,}4$.

Soit un réel $x > 1$.  Si $x\in\intervalleoo{1}{2}$ alors $\pi(x)=0 \leqslant g(x)$ et, si $x\in\intervallefo{2}{3}$ alors $\pi(x)=1 \leqslant g(\mathrm{e}) \leqslant g(x)$. Enfin, si $x\geqslant 3$, $3 \leqslant \entiere{x} \leqslant x$ donc, comme $g$ est croissante sur $\intervallefo{\mathrm{e}}{+\infty}$, $\pi(x)=\pi(\entiere{x}) \leqslant 
g(\entiere{x}) \leqslant g(x)$ ce qui achève la démonstration. \hfill $\square$

\subsubsection{Théorème de Filaseta et Trifonov}

\begin{lemme} --- Soit $d$ un entier impair. S'il existe un diviseur premier $p$ de $d$ tel que 
\[\dfrac{\ln d}{p^{v_p(d)}\ln p}+\dfrac{1}{p-1} \leqslant 1\] alors l'unique racine réelle $\xi_d$ de $\mathscr{B}_d$ est irrationnelle.
\label{lem_condition_Bessel}
\end{lemme}

\emph{Preuve}. --- Comme $\xi_d \neq0$, il est équivalent de montrer que, sous l'hypothèse du lemme, le polynôme
\[R_d:=X^d\mathscr{B}_d\left(\dfrac{2}{X}\right)=\sum_{k=0}^{d} \dfrac{(2d-k)!}{k!(d-k)!} X^k\]
n'a pas de racine rationnelle. Posons, pour tout $k\in\llbracket 0, d\rrbracket$, $c_k:=\frac{(2d-k)!}{k!(d-k)!}$ et supposons qu'il existe un diviseur premier $p$ de $d$ tel que 
\begin{equation}
\dfrac{\ln d}{p^{v_p(d)}\ln p}+\dfrac{1}{p-1} \leqslant 1.
\label{eq_Bessel}
\end{equation}  

Nous allons montrer que les coefficients $c_k$ ($k\in\llbracket 0, d\rrbracket$) de $R_d$ satisfont les hypothèses de la propriété \ref{prop_cri_irra_poly_Newton}, ce qui suffit pour conclure.

Remarquons que, pour tout $k\in\llbracket 0, d\rrbracket$, $c_k=\binom{2d-k}{d}\frac{d!}{k!}$. Ainsi, $p$ ne divise pas $c_d=1$ et, pour tout $k\in\llbracket 0, d-1\rrbracket$, $p$ divise $c_k$ car $d$ divise $\frac{d!}{k!}$.

Reste à montrer que, pour tout $k\in\llbracket 1, d\rrbracket$, $v_p(c_k)>v_p(c_0)-k$.

Soit $k\in\llbracket 1, d\rrbracket$. Alors, $c_k=c_0\times\frac{(2d-k)!}{k!(d-k)!}\times\frac{d!}{(2d)!}$ donc, d'après la propriété \ref{prop_formule_Legendre} et le corollaire \ref{coro_Legendre},
\begin{align*}
v_p(c_k)&=v_p(c_0)+v_p\left(\dfrac{d!}{(d-k)!}\right)-v_p\left(\dfrac{(2d)!}{(2d-k)!}\right)-v_p(k!) \\
&> v_p(c_0)+\sum_{j=1}^{+\infty} \left(\entiere{\dfrac{d}{p^j}} - \entiere{\dfrac{d-k}{p^j}}\right) -  \sum_{j=1}^{+\infty} \left(\entiere{\dfrac{2d}{p^j}} - \entiere{\dfrac{2d-k}{p^j}}\right) - \dfrac{k}{p-1}
\end{align*}

Posons, pour tout $j\in\N^*$ et pour tout $n\in\N$, $a_j(n)=\entiere{\frac{n}{p^j}} - \entiere{\frac{n-k}{p^j}}$. Ainsi,
\[v_p(c_k)>v_p(c_0)+\sum_{j=1}^{+\infty} (a_j(d)-a_j(2d)) - \dfrac{k}{p-1}\]
et comme, pour tout entier $j>\entiere{\frac{\ln(2d)}{\ln p}}$, $p^j>2d$, $a_j(d)=a_j(2d)=0$,
\[v_p(c_k)>v_p(c_0)+\sum_{j=1}^{\entiere{\frac{\ln(2d)}{\ln p}}} (a_j(d)-a_j(2d)) - \dfrac{k}{p-1}.\]
Pour tous entiers strictement positifs $\ell$ et $j$, $\entiere{\frac{\ell}{p^j}}$ représente le nombre de multiples de $p^j$ dans $\llbracket 1, \ell \rrbracket$ donc, pour tous entiers strictement positifs $j$ et $n$, $a_j(n)$ représente le nombre de multiples de $p^j$ dans $\llbracket n-k+1, n\rrbracket$. 

Notons $r=v_p(d)$. Remarquons que, puisque $p$ est impair, on a également $r=v_p(2d)$. 

\emph{1\up{er} cas}. --- Soit $j$ un entier tel que $1 \leqslant j \leqslant r$. \'Ecrivons $d=p^rm$ où $m$ est un entier qui n'est pas divisible par $p$. Alors,
\[a_j(d)=\entiere{\dfrac{p^rm}{p^j}}-\entiere{\dfrac{p^rm-k}{p^j}}=p^{r-j}m-\entiere{p^{r-j}m-\dfrac{k}{p^j}}=-\entiere{-\dfrac{k}{p^j}}\]
et
\[a_j(2d)=\entiere{\dfrac{2p^rm}{p^j}}-\entiere{\dfrac{2p^rm-k}{p^j}}=2p^{r-j}m-\entiere{2p^{r-j}m-\dfrac{k}{p^j}}=-\entiere{-\dfrac{k}{p^j}}\]
donc $a_j(d)-a_j(2d)=0$.

\emph{2\up{ième} cas}. --- Soit $j$ un entier tel que $j > r$. Dans $\llbracket 2d-k+1, 2d\rrbracket$, il y a $k$ entiers donc le nombre de multiples de $p^r$ est inférieur ou égal $\entiere{\frac{k}{p^r}}+1$. De plus, $2d$ est un multiple de $p^r$ mais pas de $p^j$ car $j>r$ donc $a_j(2d) \leqslant a_r(2d)-1 \leqslant \entiere{\frac{k}{p^r}}$. Il s'ensuit que 
\[a_j(d)-a_j(2d) \geqslant -a_j(2d) \geqslant -\entiere{\dfrac{k}{p^r}} \geqslant -\dfrac{k}{p^r}.\]

Il suit de ce qui précède que
\[v_p(c_k)>v_p(c_0)+\left(\sum_{j=r+1}^{\entiere{\frac{\ln(2d)}{\ln p}}} -\dfrac{k}{p^r}\right) - \dfrac{k}{p-1}=v_p(c_0)-\left(\entiere{\dfrac{\ln(2d)}{\ln p}} - r\right)\dfrac{k}{p^r}- \dfrac{k}{p-1}.\]
Or, comme $p\geqslant 3$ et $r\geqslant 1$,
\[\entiere{\dfrac{\ln(2d)}{\ln p}} - r=\entiere{\dfrac{\ln 2}{\ln p}+\dfrac{\ln d}{\ln p}} - r \leqslant 1+\entiere{\dfrac{\ln d}{\ln p}}-1 \leqslant \dfrac{\ln d}{\ln p}\]
et ainsi, grâce à \eqref{eq_Bessel},
\[v_p(c_k) > v_p(c_0)-k\left(\dfrac{\ln d}{p^r\ln p}+\dfrac{1}{p-1}\right) \geqslant v_p(c_0)-k\]
ce qui permet de conclure. \hfill $\square$

\begin{theoreme}[Filaseta \& Trifonov, 2002] --- Soit $d$ un entier impair supérieur ou égal à $3$. Alors, l'unique racine réelle $\xi_d$ de $\mathscr{B}_d$ est irrationnelle.
\end{theoreme}

\emph{Preuve}. --- On raisonne par l'absurde en supposant que $\xi_d$ est rationnelle. On note $\mathcal{P}(d)$ l'ensemble des diviseurs premiers de $d$ et, pour tout $p\in\mathcal{P}(d)$, on pose $r_p=v_p(d)$. Enfin, on note $p_m$ le maximum de $\mathcal{P}(d)$ i.e. le plus grand diviseur premier de $d$.

Comme $\xi_d$ est rationnelle, d'après le lemme \ref{lem_condition_Bessel}, pour tout $p\in\mathcal{P}(d)$, 
\begin{equation}
\dfrac{\ln d}{p^{r_p}\ln p}+\dfrac{1}{p-1} > 1.
\label{eq_Bessel_minoration}
\end{equation}
donc
\begin{equation}
p^{r_p} < \dfrac{p-1}{(p-2)\ln p} \ln d.
\label{eq_Bessel_majoration}
\end{equation}

Distinguons deux cas.

\emph{1\up{er} cas}. --- Supposons que $p_m \geqslant 29$. Alors, d'après \eqref{eq_Bessel_minoration},
\begin{equation}
\dfrac{\ln d}{p_m\ln p_m} \geqslant \dfrac{\ln d}{p_m^{r_{p_m}}\ln p_m} > 1-\dfrac{1}{p_m-1}=\dfrac{p_m-2}{p_m-1}
\label{eq_Bessel_minoration_1}
\end{equation}
donc $\ln d > \frac{p_m(p_m-2)}{p_m-1}\ln p_m$ et ainsi $d>p_m^{\frac{p_m(p_m-2)}{p_m-1}}$. Comme la fonction $x\mapsto \frac{x(x-2)}{x-1}$ est croissante sur $\intervalleoo{1}{+\infty}$ et comme $p_m\geqslant 29$, on en déduit que $d>29^{\frac{29.27}{28}}=29^{\frac{783}{28}}:=d_0$.

Montrons que $p_m \leqslant \frac{351}{250} \cdot \frac{\ln d}{\ln(\ln d)}$. Si $\ln p_m \geqslant \frac{7000}{9477}\ln(\ln d)$ alors, d'après \eqref{eq_Bessel_minoration_1}, par décroissance de $x\mapsto \frac{x-1}{x-2}$ sur $\intervalleoo{2}{+\infty}$,
\[p_m<\dfrac{p_m-1}{p_m-2}\cdot\dfrac{\ln d}{\ln p_m}< \dfrac{28}{27}\cdot\dfrac{\ln d}{\frac{7000}{9477}\ln(\ln d)}=\dfrac{351}{250}\cdot\dfrac{\ln d}{\ln(\ln d)}.\]
Sinon, $\ln p_m < \frac{7000}{9477}\ln(\ln d)$ donc $p_m<(\ln d)^{\frac{7000}{9477}}$ et il suffit de montrer que $(\ln d)^{\frac{7000}{9477}} \leqslant \frac{351}{250}\cdot\frac{\ln d}{\ln(\ln d)}$ i.e. $\frac{351}{250}\cdot\frac{(\ln d)^{\frac{2477}{9477}}}{\ln(\ln d)} \geqslant 1$. Or, $d>d_0$ donc $\ln d > \ln d_0=\frac{783}{28}\ln 29>\mathrm{e}^{\frac{9477}{2477}}$ donc, d'après le lemme \ref{lem_étude_fonction_Bessel}, $\frac{351}{250}\cdot\frac{(\ln d)^{\frac{2477}{9477}}}{\ln(\ln d)}>\frac{351}{250}\cdot\frac{(\ln d_0)^{\frac{2477}{9477}}}{\ln(\ln d_0)} > 1{,}01$ ce qui suffit pour conclure.

Par ailleurs, en posant, pour tout entier  $n\geqslant 3$, $t_n:=\ln\left(\frac{n-1}{(n-2)\ln n}\right)$, on déduit de l'inégalité \eqref{eq_Bessel_majoration} que, pour tout $p\in\mathcal{P}(d)$, $r_p\ln p < \ln(\ln d)+t_p$ donc
\[\ln d = \ln\left(\prod_{p\in \mathcal{P}(d)} p^{r_p}\right) = \sum_{p\in\mathcal{P}(d)} r_p\ln p \leqslant  \sum_{p\in\mathcal{P}(d)} (\ln(\ln d) + t_p) \leqslant \text{Card}(\mathcal{P}(d))\ln(\ln d)+\sum_{p\in\mathcal{P}(d)} t_p.\]
Or, la suite $(t_n)_{n\geqslant 3}$ est décroissante et $t_5<0$ donc, pour tout $n\geqslant 5$, $t_n<0$. Il s'ensuit que $\sum_{p\in\mathcal{P}(d)} t_p \leqslant t_3+t_{p_m} \leqslant t_3+t_{29} < 0$. D'autre part, $\text{Card}(\mathcal{P}(d)) \leqslant \pi(p_m)$ donc, comme $p_m \leqslant \frac{351}{250} \cdot \frac{\ln d}{\ln(\ln d)}$,
\[\ln d \leqslant \pi\left(\dfrac{351}{250} \cdot \dfrac{\ln d}{\ln(\ln d)}\right)\ln(\ln d).\]
Dès lors, d'après la propriété \ref{prop_majoration_pi},
\[\ln d \leqslant \dfrac{(1+\ln 4)\cdot\dfrac{351}{250} \cdot \dfrac{\ln d}{\ln(\ln d)}}{\ln\left(\dfrac{351}{250} \cdot \dfrac{\ln d}{\ln(\ln d)}\right)}\ln(\ln d) = \dfrac{(1+\ln 4)\cdot\dfrac{351}{250}}{\ln\left(\dfrac{351}{250} \cdot \dfrac{\ln d}{\ln(\ln d)}\right)}\ln d. \]
et ainsi
\[\dfrac{\ln d}{\ln(\ln d)} \leqslant \dfrac{250}{351}\mathrm{exp}\left(\dfrac{351(1+\ln 4)}{250}\right) < 20{,}4.\]
Mais, par ailleurs, comme $\ln d >\ln d_0=\frac{783}{28} \ln 29 > e$, d'après le lemme \ref{lem_étude_fonction_Bessel},
\[\dfrac{\ln d}{\ln(\ln d)}>\dfrac{\ln d_0}{\ln(\ln d_0)} > 20{,}7\]
ce qui conduit à la contradiction voulue.

\emph{2\up{ième} cas}. --- Supposons que $p_m \leqslant 23$. Alors, comme $d$ est impair, $\text{Card}(\mathcal{P}(d))\leqslant 8$ donc, d'après \eqref{eq_Bessel_majoration},
\[d=\prod_{p\in\mathcal{P}(d)} p^{r_p} \leqslant \prod_{p\in\mathcal{P}(d)} \dfrac{p-1}{(p-2)\ln p}\ln d \leqslant \prod_{p\in\mathcal{P}(d)} \dfrac{2}{\ln 3}\ln d \leqslant \left[\dfrac{2}{\ln 3}\ln d\right]^8.\]
Ainsi, $\frac{d^{\frac{1}{8}}}{\ln d} \leqslant \frac{2}{\ln 3}$. Or, $\frac{(10^{15})^{\frac{1}{8}}}{\ln(10^{15})} > \frac{2}{\ln 3}$ donc, comme $10^{15}\geqslant \mathrm{e}^8$, on déduit du lemme \ref{lem_étude_fonction_Bessel} que $d \leqslant 10^{15}$. Mais alors, d'après \eqref{eq_Bessel_majoration}, 
\begin{equation}
\dfrac{p_m(p_m-2)\ln p_m}{p_m-1} \leqslant \dfrac{p_m^{r_m}(p_m-2)\ln p_m}{p_m-1} \leqslant \ln d \leqslant 15\ln 10.
\label{eq_inégalité_pm}
\end{equation}
Or, la fonction $g:x\mapsto \frac{x(x-2)\ln x}{x-1}$ est croissante sur $\intervalleoo{1}{+\infty}$ et $g(17)>15\ln 10$ donc $p_m \leqslant 13$. Ainsi, $\text{Card}(\mathcal{P}(d))\leqslant 5$ donc, comme précédemment, $d\leqslant \left[\frac{2}{\ln 3}\ln d\right]^5$ i.e. $\frac{d^{\frac{1}{5}}}{\ln d} \leqslant \frac{2}{\ln 3}$ ce qui impose, d'après le lemme \ref{lem_étude_fonction_Bessel}, $d \leqslant 4.10^{7}$ et donc, on déduit de \eqref{eq_inégalité_pm} que $g(p_m) \leqslant \ln(4.10^7)$. Il s'ensuit, comme précédemment, que $p_m \leqslant 7$ et ainsi, grâce à  l'inégalité \eqref{eq_Bessel_majoration},
\[d \leqslant \dfrac{2}{\ln(3)}\cdot \dfrac{4}{3\ln(5)}\cdot \dfrac{6}{5\ln(7)} (\ln d)^3 \qquad \text{i.e.} \qquad \dfrac{d^{\frac{1}{3}}}{\ln d} \leqslant \left(\dfrac{16}{5\ln(3)\ln(5) \ln(7)}\right)^{\frac{1}{3}}\]
ce qui impose, d'après le lemme \ref{lem_étude_fonction_Bessel}, $d\leqslant 74$ et donc, d'après \eqref{eq_Bessel_majoration}, $g(p_m) \leqslant \ln(74)$. On en déduit que $p_m=3$ donc $d=3^{r_3}$ et ainsi, d'après \eqref{eq_Bessel_majoration}, $\frac{d}{\ln d} \leqslant \frac{2}{\ln 3} < e$ ce qui est absurde car, d'après le lemme \ref{lem_étude_fonction_Bessel}, pour tout $x>1$, $\frac{x}{\ln x} \geqslant e$.

On conclut donc que $\xi_d$ est irrationnelle. \hfill $\square$

\section{Polynômes de Bernoulli}
\subsection{Nombres de Bernoulli: définition et propriétés}

On considère la suite $(b_n)_{n\in\N}$ des nombres de Bernoulli qui sont les rationnels définis par 
\[\begin{cases} b_0=1 \\ \forall n\in\N^*,~b_n=-\dfrac{1}{n+1}\sum_{j=0}^{n-1}\dbinom{n+1}{j}b_j \end{cases}.\]
On a donc $b_1=-\frac{1}{2}$, $b_2=\frac{1}{6}$, $b_3=0$, $b_4=-\frac{1}{30}$ et $b_5=0$. 

\begin{lemme} --- La série entière $\sum \frac{b_n}{n!} x^n$ a un rayon de convergence $R\geqslant 1$ et, pour tout $x\in\intervalleoo{-R}{R}$, 
\[(\mathrm{e}^{x}-1)\sum_{n=0}^{+\infty} \frac{b_n}{n!} x^n = x.\]
\end{lemme}

\emph{Preuve}. --- Montrons par récurrence forte que, pour tout $n\in\N$, $\abs{b_n}\leqslant n!$. Comme $b_0=1$, l'inégalité est vraie pour $n=0$. On suppose que, pour un certain $n\in\N$, $\abs{b_j}\leqslant j!$ pour tout $j\in\llbracket 0, n\rrbracket$. Alors,
\begin{align*}
\abs{b_{n+1}} &\leqslant \dfrac{1}{n+2}\sum_{j=0}^{n} \dbinom{n+2}{j} \abs{b_j} \leqslant \dfrac{1}{n+2}\sum_{j=0}^{n} \dfrac{(n+2)!}{j!(n+2-j)!} j! \\
& \leqslant (n+1)!\sum_{j=0}^{n} \dfrac{1}{(n+2-j)!} \leqslant (n+1)!\sum_{k=2}^{+\infty} \dfrac{1}{k!}=(n+1)!(\mathrm{e}-2)
\end{align*}
donc, comme $\mathrm{e}-2\leqslant 1$, $\abs{b_{n+1}} \leqslant (n+1)!$ ce qui achève la récurrence.

Ainsi, la suite $\left(\frac{b_n}{n!}\right)$ est bornée donc le lemme d'Abel assure que le rayon de convergence $R$ de la série $\sum \frac{b_n}{n!}x^n$ est au moins égal à $1$.

Soit $x\in\intervalleoo{-R}{R}$. Alors, en développant la fonction exponentielle en série entière, il vient
\begin{align*}
(\mathrm{e}^x-1)\sum_{n=0}^{+\infty} \dfrac{b_n}{n!}x^n &= \sum_{n=1}^{+\infty} \dfrac{1}{n!}x^n \sum_{n=0}^{+\infty} \dfrac{b_n}{n!}x^n = \sum_{k=0}^{+\infty} \dfrac{1}{(k+1)!}x^{k+1}\sum_{j=0}^{+\infty} \dfrac{b_j}{j!}x^j \\
& = \sum_{n=0}^{+\infty} \left(\sum_{j=0}^{n} \dfrac{b_j}{j!}\times \dfrac{1}{(n-j+1)!}\right) x^{n+1} = \sum_{n=0}^{+\infty} \dfrac{1}{(n+1)!}\left(\sum_{j=0}^{n} \dbinom{n+1}{j} b_j\right) x^{n+1}.
\end{align*}
Or, pour tout $n\in\N^*$,
\[\sum_{j=0}^{n} \dbinom{n+1}{j} b_j = (n+1)b_n + \sum_{j=0}^{n-1} \dbinom{n+1}{j} b_j = (n+1)\left[b_n + \dfrac{1}{n+1}\sum_{j=0}^{n-1} \dbinom{n+1}{j} b_j\right]\]
donc, par définition,
\begin{equation}
\sum_{j=0}^{n} \dbinom{n+1}{j} b_j=(n+1)(b_n-b_n)=0
\label{eq_bernoulli}
\end{equation}
et ainsi $(\mathrm{e}^x-1)\sum\limits_{n=0}^{+\infty} \frac{b_n}{n!}x^n = \frac{1}{1!}b_0x =x$. \hfill $\square$

\begin{rem} --- On peut montrer qu'en fait $R=2\pi$ et qu'ainsi l'égalité $(\mathrm{e}^x-1)\sum\limits_{n=0}^{+\infty} \frac{b_n}{n!}x^n = x$ est vraie pour tout $x\in\intervalleoo{-2\pi}{2\pi}$. Voir, par exemple, \cite[chap. 2, § 10]{Rad73}.
\label{rem_rayon_Bernoulli}
\end{rem}

\begin{corollaire} --- Pour tout $k\in\N^*$, $b_{2k+1}=0$.
\label{coro_bernoulli_impair}
\end{corollaire}

\emph{Preuve}. --- Considérons la fonction $g : x \mapsto \sum\limits_{n=2}^{+\infty} \frac{b_n}{n!}x^n$ définie sur l'intervalle $\intervalleoo{-R}{R}$. Alors, pour tout $x\in\intervalleoo{-R}{R}\setminus\{0\}$, 
\[g(x) = \dfrac{x}{\mathrm{e}^x-1}-b_0-b_1x=\dfrac{x}{\mathrm{e}^x-1}-1+\dfrac{1}{2}x=\dfrac{x(\mathrm{e}^x+1)}{2(\mathrm{e}^x-1)}-1=\dfrac{x}{2}\coth\left(\dfrac{x}{2}\right)-1.\]
Comme la fonction coth est impaire, $x\mapsto \frac{x}{2}\coth\left(\frac{x}{2}\right)$ est paire et donc $g$ est paire. Son développement en série entière ne contient donc que des puissances paires de $x$ ce qui assure que, pour tout $n\geqslant 2$, les coefficients d'indices impairs de $(b_n)$ sont nuls i.e. pour tout $k\in\N^*$, $b_{2k+1}=0$. \hfill $\square$

\begin{propriete}[Formule de Faulhaber] --- Soit $N\in\N^*$. Pour tout $r\in\N^*$, 
\[\sum_{m=1}^{N-1} m^r = \dfrac{1}{r+1}\sum_{j=0}^{r} \dbinom{r+1}{j}b_jN^{r+1-j} = \sum_{j=0}^{r} \dfrac{1}{r+1-j}\dbinom{r}{j}b_jN^{r+1-j}.\]
\label{prop_Faulhaber}
\end{propriete}

\emph{Preuve}. --- Si $N=1$, la première égalité découle de l'identité \eqref{eq_bernoulli}.

Supposons $N\geqslant 2$. Soit $x\in\intervalleoo{-R}{R}\setminus\{0\}$. Alors, d'une part,
\[\sum_{m=0}^{N-1} \mathrm{e}^{mx} = 1+ \sum_{m=1}^{N-1} \sum_{r=0}^{+\infty} \dfrac{(mx)^{r}}{r !} = 1+ \sum_{r=0}^{+\infty}\left(\dfrac{1}{r!}\sum_{m=1}^{N-1} m^{r}\right) x^r\]
et, d'autre part, 
\[\sum_{m=0}^{N-1} \mathrm{e}^{mx} = \dfrac{\mathrm{e}^{Nx}-1}{\mathrm{e}^x-1}=\dfrac{1}{\mathrm{e}^x-1}\sum_{k=1}^{+\infty} \dfrac{(Nx)^k}{k!} = N\dfrac{x}{\mathrm{e}^x-1}\sum_{k=0}^{+\infty} \dfrac{N^k}{(k+1)!} x^k\]
donc, puisque $\frac{x}{\mathrm{e}^x-1}=\sum\limits_{n=0}^{+\infty} \frac{b_n}{n!} x^n$,
\begin{align*}
\sum_{m=0}^{N-1} \mathrm{e}^{mx} &=N \left(\sum_{n=0}^{+\infty} \dfrac{b_n}{n!} x^n\right)\left(\sum_{k=0}^{+\infty} \dfrac{N^k}{(k+1)!} x^k\right) \\
&=\sum_{r=0}^{+\infty} N \left(\sum_{j=0}^{r} \dfrac{b_j}{j!}\times \dfrac{N^{r-j}}{(r-j+1)!}\right)x^r \\
&= \sum_{r=0}^{+\infty} \left(\dfrac{1}{(r+1)!}\sum_{j=0}^{r} \dbinom{r+1}{j} b_j N^{r+1-j}\right)x^r.
\end{align*}
Ainsi, pour tout $x\in\intervalleoo{-R}{R}\setminus\{0\}$, 
\[1+\sum_{r=0}^{+\infty}\left(\dfrac{1}{r!}\sum_{m=1}^{N-1} m^{r}\right) x^r = \sum_{r=0}^{+\infty} \left(\dfrac{1}{(r+1)!}\sum_{j=0}^{r} \dbinom{r+1}{j} b_j N^{r+1-j}\right)x^r,\]
égalité qui reste vraie pour $x=0$ car $1+\frac{1}{0!}\sum\limits_{m=1}^{N-1} m^{0}= N$ et $\frac{1}{1!}\sum\limits_{j=0}^{0} \binom{1}{j} b_j N^{1-j}=b_0N=N$. On en déduit que, pour tout $r\in\N^*$,
\[\sum_{m=1}^{N-1} m^{r} = \dfrac{r!}{(r+1)!} \sum_{j=0}^{r} \dbinom{r+1}{j} b_j N^{r+1-j} = \dfrac{1}{r+1} \sum_{j=0}^{r} \dbinom{r+1}{j} b_j N^{r+1-j}.\]
Ainsi, la première égalité est démontrée pour tout $N\in\N^*$. 

La seconde égalité s'en déduit immédiatement car $\frac{1}{r+1} \binom{r+1}{j} =\frac{1}{r+1-j}\binom{r}{j}$.  \hfill $\square$

\subsection{Théorème de von Staudt et Clausen}

\begin{lemme} --- Soit $K$ un corps et $G$ un sous-groupe fini de $(K^*, \times)$. Alors, $G$ est cyclique.
\label{lem_groupeinversiblecorps}
\end{lemme}

\emph{Preuve}. --- Notons $n$ l'ordre de $G$. Alors, d'après le théorème de Lagrange, pour tout $x\in G$, $x^n=1$. Ainsi, le polynôme $P=X^n-1$ de $K[X]$ admet au moins $n$ racines distincts dans $G \subset K$. Or, ce polynôme étant de degré $n$, il admet au plus $n$ racines dans $K$. On en déduit que $G$ est exactement l'ensemble des racines de $P$ dans $K$ et que $P$ est scindé à racines simples sur $K$. 

Soit $d$ un diviseur de $n$. Il existe $m\in\N$ tel que $n=dm$ et ainsi 
\[P = \left(X^d\right)^m-1^m=(X^d-1)\sum_{i=0}^{m-1}\left(X^d\right)^i\]
donc $X^d-1$ divise $P$. Dès lors, $X^d-1$ est également scindé à racines simples et toutes ses racines sont dans $G$. Autrement dit, il y a dans $G$ exactement $d$ éléments ayant un ordre qui divise $d$. Notons $G_d$ cet ensemble.

Soit $p$ un nombre premier et $r$ un entier naturel non nul. Supposons que $p^r$ divise $n$. Alors, d'après ce qui précède, le nombre d'éléments de $G$ d'ordre $p^r$ est $\abs{G_{p^r}}-\abs{G_{p^{r-1}}}=p^r-p^{r-1}>0$. Ainsi, il existe dans $G$ au moins un élément d'ordre $p^r$. \'Ecrivons alors la décomposition de $n$ en produit de facteurs premiers: $n=\prod\limits_{i=1}^{j} p_i^{r_i}$. Pour tout $i\in\llbracket1, j\rrbracket$, il existe un élément $x_i\in G$ d'ordre $p_i^{r_i}$. Posons $x:=\prod\limits_{i=1}^j x_i \in G$ et notons $q$ l'ordre de $x$. D'une part, comme $x\in G$, $q$ divise $n$. D'autre part, considérons $k\in\llbracket 1, j\rrbracket$ et notons $n_k=\frac{n}{p_k^{r_k}}\in\N$. Alors, comme $G$ est un groupe commutatif (car $K$ est un corps),
\[1=(x^q)^{n_k}=\left(\prod_{i=1}^j x_i\right)^{qn_k}=\prod_{i=1}^j (x_i^{n_k})^{q}=x_k^{qn_k}\]
car, pour tout $i\in\llbracket 1, j\rrbracket$ tel que $i\neq k$, $x_i^{n_k}=1$ puisque $p_i^{r_i}$ divise $n_k$. Dès lors, $p_k^{r_k}$ divise $qn_k$ et, comme $p_k^{r_k}$ est premier avec $n_k$, par le lemme de Gauss, $p_k^{r_k}$ divise $q$. Comme les entiers $p_i^{r_i}$ pour $i\in\llbracket 1, j \rrbracket$ sont premiers entre eux deux à deux, on conclut que $n=\prod\limits_{i=1}^{j} p_i^{r_i}$ divise $q$ et, finalement, $n=q$.  Ainsi, $G$ admet un élément d'ordre $n$: il est donc cyclique. \hfill $\square$

\bigskip

\begin{lemme} --- On définit, pour tout $r\in\N^*$, la fonction $\varepsilon_{r}$ sur $\mathbb{P}$  par $\varepsilon_{r}(p) = \begin{cases} 1\text{ si $p-1$ divise $r$} \\ 0\text{ sinon}. \end{cases}$

Alors, pour tout $k\in\N^*$ et tout $p\in\mathbb{P}$, $\sum\limits_{m=1}^{p-1} m^{r} \equiv -\varepsilon_{r}(p) \pmod{p}$.
\label{lem_prop_epsilon_p}
\end{lemme}

\emph{Preuve}. --- Soit $r\in\N^*$ et $p\in\mathbb{P}$. Si $p-1$ divise $r$ alors il existe un entier naturel $\ell$ tel que $r=\ell(p-1)$. Or, d'après le petit théorème de Fermat, pour tout $m\in\llbracket 1, p-1\rrbracket$, $m^{p-1}\equiv 1 \pmod{p}$ donc $m^{r}=\left(m^{p-1}\right)^{\ell} \equiv 1 \pmod{p}$. Ainsi, $\sum\limits_{m=1}^{p-1} m^{r} \equiv p-1 \pmod{p} \equiv -1 \pmod{p} \equiv -\varepsilon_{r}(p) \pmod{p}$. Supposons à présent que $p-1$ ne divise pas $r$.  Comme $p$ est premier, $\EQ{\Z}{p\Z}$ est un corps donc, d'après le lemme \ref{lem_groupeinversiblecorps}, le groupe $(\EQ{\Z}{p\Z})^*$ est cyclique. Il existe donc un entier $a\in\llbracket 1, p-1\rrbracket$ tel que la classe $\overline{a}$ de $a$ dans $(\EQ{\Z}{p\Z})^*$ engendre $(\EQ{\Z}{p\Z})^*$. L'ordre de $\overline{a}$ est alors $p-1$ donc $\overline{a}^{r} \neq \overline{1}$ car sinon $p-1$ diviserait $r$. De plus, l'application $\overline{x} \mapsto \overline{ax}$ est une bijection de $(\EQ{\Z}{p\Z})^*$ dans lui-même donc  
\[ \overline{a}^{r} \sum_{m=1}^{p-1} \overline{m}^{r}=\sum_{m=1}^{p-1} (\overline{am})^{r} = \sum_{m=1}^{p-1} \overline{m}^{r}.\]
Ainsi, $(\overline{a}^{r}-\overline{1})\sum\limits_{m=1}^{p-1} \overline{m}^{r} = \overline{0}$ et $\overline{a}^{r}-\overline{1}\neq \overline{0}$ donc, par intégrité de $\EQ{\Z}{p\Z}$, $\sum\limits_{m=1}^{p-1} \overline{m}^{r} = \overline{0}$, ce qui revient à dire que $\sum\limits_{m=1}^{p-1} m^{r}\equiv 0 \pmod{p} \equiv -\varepsilon_{r}(p) \pmod{p}$.

Ainsi, l'égalité est établie dans tous les cas. \hfill $\square$

\begin{theoreme}[von Staudt et Clausen, 1840] --- Soit $n\in\N^*$. Si $n=1$ ou si $n$ est pair alors $b_{n}+\sum\limits_{p\in\mathbb{P}} \frac{\varepsilon_{n}(p)}{p}$ est un entier.
\label{theo_vS_C}
\end{theoreme}

\emph{Preuve}. --- Remarquons que la somme est évidemment finie puisque $n$ n'a qu'un nombre fini de diviseurs donc $\varepsilon_n(p)\neq 0$ pour un nombre fini de nombres premiers $p$.

Montrons par récurrence forte que, pour tout $n\in\N^*$, si $n=1$ ou $n$ est pair alors $b_{n}+\sum\limits_{p\in\mathbb{P}} \frac{\varepsilon_n(p)}{p}$ est un entier.

Si $n=1$ alors, pour tout nombre premier $p\geqslant 3$, $\varepsilon_n(p)=0$ donc $b_{n}+\sum\limits_{p\in\mathbb{P}} \frac{\varepsilon_n(p)}{p} = b_1+\frac{1}{2}=-\frac{1}{2}+\frac{1}{2}=0$ est bien un entier. De même, si $n=2$ alors, pour tout nombre premier $p\geqslant 5$, $\varepsilon_n(p)=0$ donc $b_{n}+\sum\limits_{p\in\mathbb{P}} \frac{\varepsilon_n(p)}{p} = b_2+\frac{1}{2}+\frac{1}{3}=\frac{1}{6}+\frac{1}{2}+\frac{1}{3}=1$ est bien un entier.

Supposons que, pour un certain entier $k\geqslant 1$, $b_{j}+\sum\limits_{p\in\mathbb{P}} \frac{\varepsilon_j(p)}{p}$ soit un nombre entier pour tout nombre $j\in\{1,2,4,6,..., 2k\}$. Soit $q\in\mathbb{P}$. En appliquant la propriété \ref{prop_Faulhaber} et le lemme \ref{lem_prop_epsilon_p} avec $N=q$ et $r=2(k+1)$, il vient  
\[\varepsilon_{2(k+1)}(q) + \sum_{j=0}^{2(k+1)} \dfrac{1}{2(k+1)+1-j}\dbinom{2(k+1)}{j} b_jq^{2(k+1)+1-j} = \varepsilon_{2(k+1)}(q)+\sum_{m=1}^{q-1} m^{2(k+1)} \equiv 0 \pmod{q}.\]
Ainsi, en divisant par $q$ et en tenant compte du fait que $b_{2k+1}=0$ d'après le corollaire \ref{coro_bernoulli_impair}, on est assuré que le nombre
\[T_q:=b_{2k+2} + \dfrac{\varepsilon_{2k+2}(q)}{q} + \sum_{j=0}^{2k} \dfrac{1}{2k+3-j}\dbinom{2k+2}{j} b_jq^{2k+2-j}\]
 est un entier.
 
Soit $j\in\llbracket 0, 2k\rrbracket$. \'Ecrivons le nombre rationnel $\frac{1}{2k+3-j}\binom{2k+2}{j} b_jq^{2k+2-j}$ sous forme irréductible $\frac{M_j}{N_j}$. On va montrer que $q$ ne divise pas $N_j$.

Si $j$ est un nombre impair différent de $1$, $b_j=0$ donc $N_j=1$ n'est pas divisible par $q$. Sinon, $j\in\{1, 2, 4, ..., 2k\}$ donc, par hypothèse de récurrence, il existe un entier $K_j$ tel que $b_j=K_j-\sum\limits_{p\in\mathbb{P}} \frac{\varepsilon_j(p)}{p}$ et ainsi $qb_j=qK_j-\varepsilon_j(q)-\sum\limits_{p\in\mathbb{P}\setminus\{q\}}\frac{q\varepsilon_j(p)}{p}$ donc $qb_j$ s'écrit sous forme irréductible $\frac{\gamma_j}{\delta_j}$ avec $q \nmid \delta_j$. \'Etant donné que
\[\dfrac{M_j}{N_j}=\dfrac{1}{2k+3-j}\dbinom{2k+2}{j} (qb_j)q^{2k+1-j}= \dfrac{1}{2k+3-j}\dbinom{2k+2}{j} \dfrac{\gamma_j}{\delta_j} q^{2k+1-j},\]
on a, en utilisant le lemme \ref{lem_majoration_vp} avec $d=2k+3-j \geqslant 3$,
\[v_q\left(N_j\right) \leqslant v_q(2k+3-j)-(2k+1-j) \leqslant 2k+1-j-(2k+1-j)=0\]
donc $q$ ne divise pas $N_j$.
 
On conclut donc que, dans tous les cas, $q$ ne divise par $N_j$. 
 
 On a donc montré qu'il existe un entier $T_q$ tel que
 \[b_{2k+2} + \dfrac{\varepsilon_{2k+2}(q)}{q}=T_q-\sum_{j=0}^{2k} \dfrac{M_j}{N_j}\]
 où $M_j$ et $N_j$ sont des entiers tels que $q$ ne divise pas $N_j$ donc le dénominateur de $b_{2k+2} + \frac{\varepsilon_{2k+2}(q)}{q}$ n'est pas divisible par $q$ et ainsi le dénominateur de $b_{2k+2} + \sum\limits_{p\in\mathbb{P}} \frac{\varepsilon_{2k+2}(p)}{p}$ n'est pas divisible par $q$. Comme ceci est vrai pour tout nombre premier $q$, on conclut que $b_{2k+2} + \sum\limits_{p\in\mathbb{P}} \frac{\varepsilon_{2k+2}(p)}{p}$ est un entier. Le résultat est donc bien montré par récurrence. \hfill $\square$

\begin{rem} --- Le théorème précédent a été démontré indépendamment par Clausen \cite{Cla40} et von Staudt \cite{vSt40} en 1840. La démonstration que nous en avons donné est due à Rado \cite{Rad34}.
\label{rem_preuve_vSC}
\end{rem}

\begin{corollaire} --- Pour tout $n\in\N$, on note $\mathcal{P}_n$ l'ensemble des nombres premiers tels que $p-1$ divise $n$ et $\beta_n$ le dénominateur de $b_n$. Si $n=1$ ou $n$ est un entier pair strictement positif alors $\beta_n=\prod\limits_{p\in\mathcal{P}_n} p$ et si $n=0$ ou $n$ est un entier impair au moins égal à 3 alors $\beta_n=1$.
\label{coro_vS_C}
\end{corollaire}

\emph{Preuve}. --- Comme $b_1=-\frac{1}{2}$, $\beta_1=2=\prod\limits_{p\in\mathcal{P}_1} p$ car $\mathcal{P}_1=\{2\}$. Soit $n$ est un entier pair strictement positif. Alors, d'après le théorème \ref{theo_vS_C}, il existe un entier $E_n$ tel que $b_n=E_n-\sum\limits_{p\in\mathcal{P}_n} \frac{1}{p}$. Or, le dénominateur de $\sum\limits_{p\in\mathcal{P}_n} \frac{1}{p}$ est le P.P.C.M. des éléments de $\mathcal{P}_n$ qui n'est autre que $\prod\limits_{p\in\mathcal{P}_n} p$ car les éléments de $\mathcal{P}_n$ sont premiers entre eux deux à deux. Par suite, on a également $\beta_n=\prod\limits_{p\in\mathcal{P}_n} p$. 

Par ailleurs, $b_0=1$ donc $\beta_0=1$ et si $n$ est un entier impair au moins égal à $3$ alors, d'après le corollaire \ref{coro_bernoulli_impair}, $b_n=0$ donc $\beta_n=1$.  \hfill $\square$

\subsection{Polynômes de Bernoulli: définition et propriétés}

\begin{definition} --- On définit, pour tout $d\in\N$, le polynôme de Bernoulli d'indice $d$ par
\[B_d:=\sum_{k=0}^d \dbinom{d}{k} b_{d-k} X^{k}.\]
\end{definition}

Les cinq premiers polynômes de Bernoulli sont $B_0=1$, $B_1=X-\frac{1}{2}$, $B_2=X^2-X+\dfrac{1}{6}$, $B_3=X^3-\frac{3}{2}X^2+\frac{1}{2}X$ et $B_4=X^4-2X^3+X^2-\frac{1}{30}$.

\begin{propriete} --- Pour tout $d\in\N$, $B_d(X+1)-B_d=dX^{d-1}$ et, pour tout entier $d\geqslant 2$, $B_d(0)=B_d(1)$.
\label{prop_différence_Bernoulli}
\end{propriete}

\emph{Preuve}. --- Soit $d\in\N$. Si $d=0$ alors $B_d(X+1)-B_d=1-1=0=dX^{d-1}$. 

Supposons à présent $d\geqslant 1$.
\begin{align*} 
B_{d}(X+1)-B_d &= \sum_{k=1}^{d} \dbinom{d}{k}b_{d-k}\left[(X+1)^k-X^k\right] \\
& =  \sum_{k=1}^{d} \dbinom{d}{k}b_{d-k} \left[\sum_{j=0}^{k-1} \dbinom{k}{j} X^j\right] \\
& = \sum_{k=1}^{d}\sum_{j=0}^{k-1} \dbinom{d}{k}\dbinom{k}{j} b_{d-k}   X^j.
\end{align*}
Or, lorsque $\couple{k}{j}$ décrit $\llbracket 1, d\rrbracket \times \llbracket 0, k-1\rrbracket$, $\couple{j}{k}$ décrit $\llbracket 0, d-1\rrbracket \times \llbracket j+1, d\rrbracket$ donc
\begin{align*}
B_{d}(X+1)-B_d &=\sum_{j=0}^{d-1}\left[\sum_{k=j+1}^{d} \dfrac{d!}{k!(d-k)!}\dfrac{k!}{j!(k-j)!} b_{d-k}\right]X^j \\
&=\sum_{j=0}^{d-1}\left[\sum_{k=j+1}^{d} \dfrac{d!}{j!(d-j)!}\dfrac{(d-j)!}{(d-k)!((d-j)-(d-k))!} b_{d-k}\right]X^j \\
&=\sum_{j=0}^{d-1} \dbinom{d}{j}\left[\sum_{k=j+1}^{d} \dbinom{d-j}{d-k} b_{d-k}\right]X^j.
\end{align*}
Or, pour tout $j\in\llbracket 0, d-1\rrbracket$, en posant $\ell=d-k$,
\begin{align*}
\sum_{k=j+1}^{d} \dbinom{d-j}{d-k} b_{d-k}=\sum_{\ell=0}^{d-j-1} \dbinom{d-j}{\ell} b_{\ell}
\end{align*}
et, d'après \eqref{eq_bernoulli}, cette somme est nulle si $d-j-1\geqslant 1$ i.e. si $j\leqslant d-2$. Ainsi, le seul monôme non nul dans $B_d(X+1)-B_d$ est celui de degré $j=d-1$ et
\[B_{d}(X+1)-B_d=\dbinom{d}{d-1}\dbinom{1}{0}b_0X^{d-1}=dX^{d-1}.\]

Supposons à présent que $d\geqslant 2$. Alors, $B_d(1)-B_d(0)=B_d(0+1)-B_d(0)=d\times 0^{d-1}=0$ donc $B_d(0)=B_d(1)$. \hfill $\square$

\begin{propriete} --- Pour tout $d\in\N^*$, 
\[B_d'=dB_{d-1} \qquad \text{ et } \qquad \int_{0}^1 B_d(t) \mathrm{d}t = 0.\]
\label{prop_dérivée_Bernoulli}
\end{propriete}

\emph{Preuve}. --- Soit $d\in\N^*$. Alors,
\[B_d'=\sum_{k=1}^d k\dbinom{d}{k} b_{d-k} X^{k-1}=\sum_{k=1}^d d\dbinom{d-1}{k-1} b_{d-k} X^{k-1} = d\sum_{k=0}^{d-1} \dbinom{d-1}{k} b_{d-1-k} X^{k}\]
donc $B_d'=dB_{d-1}$.

Dès lors,
\[\int_0^1 B_d(t) \mathrm{d}t = \left[\dfrac{1}{d+1}B_{d+1}(t)\right]_0^1=\dfrac{1}{d+1}\left[B_{d+1}(1)-B_{d+1}(0)\right].\]
Or, comme $d\geqslant 1$, $d+1\geqslant 2$ donc, d'après la propriété \ref{prop_différence_Bernoulli}, $B_{d+1}(0)=B_{d+1}(1)$ et on conclut donc que 
\[\int_0^1 B_d(t) \mathrm{d}t = 0\]
comme annoncé. \hfill $\square$

\begin{propriete} --- Pour tout $d\in\N$,
\begin{enumerate}
\item $B_d(1-X)=(-1)^dB_d$;
\item $B_d=2^{d-1}\left[B_d\left(\frac{X}{2}\right)+B_d\left(\frac{X+1}{2}\right)\right]$ et $B_d\left(\frac{1}{2}\right)=\frac{1-2^{d-1}}{2^{d-1}}b_d$.
\end{enumerate}
\label{prop_Bernoulli_un_demi}
\end{propriete}

\emph{Preuve} \hspace{1cm}
\begin{enumerate}
\item On raisonne par récurrence sur $d$. L'égalité est évidente pour $d\in\{0, 1\}$ car $B_0=1$ et $B_1=X-\frac{1}{2}$. Supposons  l'égalité vraie pour un certain $d\in\N^*$. Alors,
\[B_{d+1}(1-X)'=-B_{d+1}'(1-X)=-(d+1)B_{d}(1-X)=-(d+1)(-1)^{d}B_{d}=(-1)^{d+1}B_{d+1}'\]
donc il existe une constante réelle $c$ telle que $B_{d+1}(1-X)=(-1)^{d+1}B_{d+1}+c$. En particulier, $B_{d+1}(1)=(-1)^{d+1}B_{d+1}(0)+c$. Or, comme $d+1\geq 2$, la propriété \ref{prop_différence_Bernoulli} assure que $B_{d+1}(0)=B_{d+1}(1)$ et ainsi $c=[1+(-1)^d)]B_{d+1}(0)=[1+(-1)^d]b_{d+1}$. Or, si $d$ est pair, d'après le corollaire \ref{coro_bernoulli_impair}, $b_{d+1}=0$ et, si $d$ est impair, $1+(-1)^d=0$ donc, dans tous les cas, $c=0$ et ainsi $B_{d+1}(1-X)=(-1)^{d+1}B_{d+1}$ ce qui achève la récurrence.
\item Posons, pour tout $d\in\N$, $P_d=2^{d-1}\left[B_d\left(\frac{X}{2}\right) + B_d\left(\frac{X+1}{2}\right)\right]$. Montrons par récurrence que, pour tout $d\in\N$, $B_d=P_d$. On a $P_0=2^{-1}(1+1)=1=B_0$ donc l'égalité est vraie pour $d=0$. Supposons que, pour un certain $d\in\N$, $B_d=P_d$. Alors, d'après la propriété \ref{prop_dérivée_Bernoulli},
\begin{align*}
P_{d+1}' &=2^{d}\left[\dfrac{1}{2}B_{d+1}'\left(\dfrac{X}{2}\right) + \dfrac{1}{2}B_{d+1}'\left(\dfrac{X+1}{2}\right)\right] = 2^{d-1}\left[(d+1)B_{d}\left(\dfrac{X}{2}\right) +(d+1)B_{d}\left(\dfrac{X+1}{2}\right)\right] \\
&=(d+1)P_d=(d+1)B_d=B_{d+1}'
\end{align*}
Il existe donc une constante réelle $c$ telle que $P_{d+1}=B_{d+1}+c$. Or,
\begin{align*}
\int_{0}^1 P_{d+1}(t)\mathrm{d}t &= 2^{d}\int_{0}^1 B_{d+1}\left(\dfrac{t}{2}\right) + B_{d+1}\left(\dfrac{t+1}{2}\right) \mathrm{d}t = \dfrac{2^{d+1}}{d+2}\left[B_{d+2}\left(\dfrac{t}{2}\right) + B_{d+2}\left(\dfrac{t+1}{2}\right)\right]_0^1 \\
&=\dfrac{2^{d+1}}{d+2}\left[B_{d+2}\left(\dfrac{1}{2}\right) + B_{d+2}(1) - B_{d+2}(0) - B_{d+2}\left(\dfrac{1}{2}\right) \right] = 0
\end{align*}
car $d+2\geqslant 2$ donc $B_{d+2}(0)=B_{d+2}(1)$ d'après la propriété \ref{prop_différence_Bernoulli}.

Mais, par ailleurs, 
\[\int_{0}^1 B_{d+1}(t) + c~\mathrm{d}t = c+ \int_{0}^1 B_{d+1}(t) \mathrm{d}t = c\]
d'après la propriété \ref{prop_dérivée_Bernoulli}. On conclut donc que $c=0$ i.e. $P_{d+1}=B_{d+1}$ ce qui achève la récurrence.

Soit $d\in\N$. On déduit de ce qui précède que $B_d(0)=2^{d-1}\left[B_d(0)+B_d\left(\frac{1}{2}\right)\right]$ donc $B_d\left(\frac{1}{2}\right)=\left(\frac{1}{2^{d-1}}-1\right)B(0)$ i.e. $B_d\left(\frac{1}{2}\right)=\frac{1-2^{d-1}}{2^{d-1}}b_d$. \hfill $\square$
\end{enumerate}

\subsection{Irrationalité des racines: le théorème d'Inkeri}

\begin{lemme} --- Soit un entier $d\geqslant 2$. Si $\frac{s}{t}$ est une racine rationnelle de $B_d$ écrite sous forme irréductible alors $t\in\{1, 2\}$.
\label{lem_Inkeri}
\end{lemme}

\emph{Preuve}. --- Soit $\frac{s}{t}$ un rationnel écrit sous forme irréductible et tel que $\sum\limits_{k=0}^d \binom{d}{k} b_{d-k}\left(\frac{s}{t}\right)^k=0$. Alors, en multipliant par $t^{d}$, il vient 
\[s^d+\sum_{k=0}^{d-2} \dbinom{d}{k} b_{d-k} s^k t^{d-k} = -b_1ds^{d-1}t = -\dfrac{t}{2}ds^{d-1}.\]
Supposons $t\geqslant 2$ et considérons un diviseur premier $p$ de $t$. Alors, comme $s$ et $t$ sont premiers entre eux, $p$ ne divise pas $s^d$. De plus, d'après le corollaire \ref{coro_vS_C}, pour tout $n\in\N$, si $\beta_n$ est le dénominateur de $b_n$ alors $v_p(\beta_n)\leqslant 1$ et, pour tout $k\in\llbracket 0, d-2\rrbracket$, $v_p(t^{d-k})\geqslant d-k \geqslant 2$ donc il existe deux entiers $M_k$ et $N_k$ telles que $\binom{d}{k} b_{d-k} s^k t^{d-k}=\frac{M_k}{N_k}$ avec $p \mid M_k$ et $p \nmid N_k$. On en déduit, en notant $D=\prod\limits_{k=0}^{d-2} N_k$ et en posant, pour tout $k\in\llbracket 0, d-2\rrbracket$, $C_k=\frac{D}{N_k}\in\Z$, que
\[s^dD+\sum_{k=0}^{d-2} M_kC_k=-\dfrac{t}{2}ds^{d-1}D.\]
Pour tout $k\in\llbracket 0, d-2\rrbracket$, $p \mid M_k$ donc $\sum\limits_{k=0}^{d-2} M_kC_k$ est un entier divisible par $p$. De plus, ni $s^d$ ni $D$ n'est divisible par $p$, donc $s^dD$ n'est pas divisible par $p$ et ainsi $s^dD+\sum\limits_{k=0}^{d-2} M_kC_k$ est un entier premier avec $p$. Il s'ensuit que $v_p\left(-\frac{t}{2}ds^{d-1}M\right)=0$. Par suite, $1 \leqslant v_p(t) \leqslant v_p(2 )\leqslant 1$ donc $p=t=2$.

Finalement, on conclut que $t\in\{1, 2\}$. 

\begin{theoreme}[Inkeri, 1959] --- Soit un entier $d\geqslant 2$. Le polynôme $B_d$ n'a pas d'autres racines rationnelles que $0$, $\frac{1}{2}$ et $1$ si $d$ est impair et n'a aucune racine rationnelle si $d$ est pair.
\end{theoreme}

\emph{Preuve}. --- D'après le lemme \ref{lem_Inkeri}, si $r$ est une racine rationnelle de $B_d$ alors il existe un entier $m\in\Z$ tel que $r=m$ ou $r=\frac{1}{2}+m$. De plus, d'après le point 1. de la propriété \ref{prop_Bernoulli_un_demi}, pour tout réel $r$, $B_d\left(1-r\right)=(-1)^dB_d\left(r\right)$ donc on peut se restreindre à étudier les rationnels de la forme $r=m$ ou $r=\frac{1}{2}+m$ avec $m\in\N$.

D'après la propriété \ref{prop_différence_Bernoulli}, pour tout $x\in\R$ et tout $m\in\N$,
\begin{equation}
B_d(x+m)=B_d(x) + \sum_{j=0}^{m-1} \left[B_d(x+j+1)-B_d(x+j)\right] = B_d(x)+\sum_{j=0}^{m-1} d(x+j)^{d-1}.
\label{eq_bernoulli_somme_1}
\end{equation}
Supposons que $d$ est pair. Notons, comme précédemment, $\mathcal{P}_d$ l'ensemble des nombres premiers $p$ tels que $p-1$ divise $d$. Alors, comme $d$ est pair, $3$ appartient à $\mathcal{P}_d$ donc, d'après le corollaire \ref{coro_vS_C}, le dénominateur $\beta_d$ de $b_d$ est divisible par $3$. En particulier, $b_d$ n'est pas un entier. Or, en appliquant \eqref{eq_bernoulli_somme_1} avec $x=0$, il vient, pour tout $m\in\N$, 
\[B_d(m)=B_d(0)+d\sum_{j=0}^{m-1} j^{d-1}=b_d+d\sum_{j=0}^{m-1} j^{d-1}\]
et donc, comme $d\sum\limits_{j=0}^{m-1} j^{d-1}$ est entier, on conclut que $B_d(m)$ n'est pas entier. En particulier, pour tout $m\in\N$, $B_d(m)\neq 0$. De plus, en appliquant \eqref{eq_bernoulli_somme_1} avec $x=\frac{1}{2}$, on obtient grâce au point 2. de la propriété \ref{prop_Bernoulli_un_demi}, pour tout $m\in\N$,
\[B_d\left(\dfrac{1}{2}+m\right)=B_d\left(\dfrac{1}{2}\right)+d\sum_{j=0}^{m-1} \left(\dfrac{1}{2}+j\right)^{d-1}=\dfrac{1-2^{d-1}}{2^{d-1}}b_{d}+d\sum_{j=0}^{m-1} \left(\dfrac{1}{2}+j\right)^{d-1}.\] 
Notons que, comme $d-1$ est impair, $1-2^{d-1}\equiv1-(-1)^{d-1} \pmod{3} \equiv 2 \pmod{3}$ donc $3$ ne divise pas $1-2^{d-1}$. Ainsi, $3$ divise le dénominateur de $\frac{1-2^{d-1}}{2^{d-1}}b_{d}$. Par ailleurs, le dénominateur de $d\sum\limits_{j=0}^{m-1} \left(\frac{1}{2}+j\right)^{d-1}$ est une puissance de $2$ donc $3$ divise le dénominateur de $B_d\left(m+\frac{1}{2}\right)$. En particulier, $B_d\left(\frac{1}{2}+m\right) \neq 0$. On conclut donc que, si $d$ est pair, $B_d$ n'a pas de racine rationnelle.

Supposons à présent $d$ impair. Alors, comme $d\geqslant 3$, d'après la propriété \ref{prop_différence_Bernoulli} et le corollaire \ref{coro_bernoulli_impair}, $B_{d}(1)=B_d(0)=b_d=0$. De plus, d'après le point 2. de la propriété \ref{prop_Bernoulli_un_demi}, $B_d\left(\frac{1}{2}\right)=\frac{1-2^{d-1}}{2^{d-1}}b_{d}=0$ car $b_d=0$. Ainsi, $0$, $\frac{1}{2}$ et $1$ sont bien des racines de $B_d$. Remarquons que, d'après la propriété \ref{prop_différence_Bernoulli}, $B_d\left(\frac{3}{2}\right)=B_d\left(\frac{1}{2}\right)+d\left(\frac{1}{2}\right)^{d-1}=\frac{d}{2^{d-1}}\neq 0$. Notons, ensuite, $\alpha$ l'un des deux nombres $0$ ou $\frac{1}{2}$ et considérons un entier $m \geqslant 2$. Alors, en utilisant \eqref{eq_bernoulli_somme_1} avec $x=\alpha$, on obtient 
\[B_d(\alpha+m)=B_d(\alpha)+d\sum_{j=0}^{m-1} (\alpha+j)^{d-1}=d\sum_{j=0}^{m-1} (\alpha+j)^{d-1} \geqslant d(\alpha+1)^{d-1}>0\]
donc $B_d(m+\alpha) \neq 0$. Ainsi, exceptés $0$, $\frac{1}{2}$ et $1$, aucun des nombres de la forme $m$ ou $\frac{1}{2}+m$ avec $m\in\N$ n'est racine de $B_d$. On conclut donc que, si $d$ est impair, les seules racines rationnelles de $B_d$ sont $0$, $\frac{1}{2}$ et $1$. \hfill $\square$

\begin{rem}. --- La connaissance des racines réelles des polynômes de Bernoulli, tout comme leur éventuelle irréductibilité sur $\Q$, reste partielle. 

Il n'est pas difficile de voir en utilisant la propriété \ref{prop_dérivée_Bernoulli} que, pour tout entier $k\geqslant 1$, $B_{2k+1}$ admet exactement trois racines simples dans $\intervalleff{0}{1}$, à savoir $0$, $\frac{1}{2}$ et $1$ et $B_{2k}$ admet exactement deux racines simples dans $\intervalleff{0}{1}$, $\alpha_k<\beta_k$ (avec $\beta_k=1-\alpha_k$ d'après le point 1. de la propriété \ref{prop_Bernoulli_un_demi}). Lense \cite{Len34} a montré que la suite $(\alpha_k)$ est croissante et converge vers $\frac{1}{4}$ et Lehmer \cite{Leh40} a précisé ceci en donnant un développement asymptotique de $\alpha_k$. Pour établir leurs résultats, ces deux auteurs utilisent des développements en série de Fourier qui ne sont valables que sur l'intervalle $\intervalleff{0}{1}$. Dans son article \cite{Ink59}, Inkeri a donné une autre méthode pour retrouver ces résultats sans l'aide des séries de Fourier et il a pu appliquer celle-ci aux racines strictement supérieures à $1$. C'est également dans cet article qu'il a démontré le théorème précédent sur l'irrationalité des racines des polynômes de Bernoulli. Brillhart, quand à lui, a montré dans \cite{Bri69} que, pour tout entier naturel $k$, les racines de $B_{2k+1}$ sont simples et celles de $B_{2k}$ sont de multiplicité au plus 2. On ne connait cependant pas de polynômes de Bernoulli ayant des racines multiples. Sur ces sujets, on pourra également consulter \cite{Del91} et \cite{Dil88}. Notons enfin que, récemment, Edwards et Lemming \cite{EL02} ont donné un algorithme permettant de déterminer le nombre de racines réelles d'un polynôme de Bernoulli de degré donné.

Pour ce qui est de l'irréductibilité, Carlitz \cite{Car52} et Brillhart \cite{Bri69} a montré que, sous certaines hypothèses vérifiées par l'entier $k$, le polynôme $\frac{B_{2k+1}}{X(X-1)\left(X-\frac{1}{2}\right)}$ est irréductible sur $\Q$. On conjecture que ceci est vrai pour tout $k\geq 3$ et $k\neq 11$ car Brillhart a également remarqué que $X^2-X+1$ divise $B_{11}$. De même, on conjecture que, pour tout $k\geqslant 1$, $B_{2k}$ est irréductible sur $\Q$, ce qui a été démontré, dans certains cas particuliers seulement, par McCarthy \cite{McC61}. 
\label{rem_Bernoulli_irré}
\end{rem}

\section{Polynômes d'Euler}
\subsection{Définition et propriétés}

\begin{definition} --- On définit, pour tout $d\in\N$, le polynôme d'Euler d'indice $d$ par
\[E_d:=\sum_{k=0}^d \dbinom{d}{k} \dfrac{2(1-2^{d+1-k})}{d+1-k}b_{d+1-k} X^{k}\]
où, pour tout $k\in\llbracket 0, d\rrbracket$, $b_{d-k+1}$ est le nombre de Bernoulli d'indice $d-k+1$.
\end{definition}

Il est clair que, pour tout $d\in\N$, $E_d$ est un polynôme unitaire de degré $d$ de $\Q[X]$ (car $b_1=-\frac{1}{2}$). Les cinq premiers polynômes d'Euler sont $E_0=1$, $E_1=X-\frac{1}{2}$, $E_2=X^2-X$, $E_3=X^3-\frac{3}{2}X^2+\frac{1}{4}$, $E_4=X^4-2X^3+X$.

\bigskip

La définition que nous avons choisie des polynômes d'Euler fait intervenir les nombres de Bernoulli. On peut en fait établir une relation explicite entre ces polynômes et ceux de Bernoulli comme le montre la proposition suivante. 

\begin{propriete} --- Pour tout $d\in\N$, $E_d=\frac{2}{d+1}\left[B_{d+1}-2^{d+1}B_{d+1}\left(\frac{X}{2}\right)\right]$.
\label{prop_lien_Euler_Bernoulli}
\end{propriete}

\emph{Preuve}. --- Soit $d\in\N$. Alors, comme $1-2^{d-k+1}$ s'annule pour $k=d+1$ et comme, pour tout $k\in\llbracket 0, d\rrbracket$, $\frac{1}{d+1-k}\binom{d}{k}=\frac{1}{d+1}\binom{d+1}{k}$, il vient
\begin{align*}
E_d&=\dfrac{2}{d+1}\sum_{k=0}^{d+1} \dbinom{d+1}{k} (1-2^{d+1-k})b_{d+1-k} X^{k} \\
& = \dfrac{2}{d+1}\left[\sum_{k=0}^{d+1} \dbinom{d+1}{k} b_{d+1-k} X^{k} -2^{d+1}\sum_{k=0}^{d+1} \dbinom{d+1}{k} b_{d-k+1} \left(\dfrac{X}{2}\right)^{k}\right] \\
&=\dfrac{2}{d+1}\left[B_{d+1}-2^{d+1}B_{d+1}\left(\dfrac{X}{2}\right)\right].
\end{align*}
\hfill $\square$

On en déduit les propriétés suivantes des polynômes d'Euler.

\begin{propriete} \hspace{1cm}

\begin{enumerate}
\item Pour tout $d\in\N$, $E_d(X+1)+E_d=2X^{d}$.
\item Pour tout $d\in\N^*$, $E_d(1)=-E_d(0)$. 
\item Pour tout $d\in\N^*$, $E_d'=dE_{d-1}$.
\item Pour tout $d\in\N$ et tout $a\in\R$, $E_d(X)=\sum\limits_{k=0}^d \binom{d}{k}E_{d-k}(a)(X-a)^k$.
\end{enumerate}
\label{prop_polynômes_Euler}
\end{propriete}

\emph{Preuve} \hspace{1cm}

\begin{enumerate}
\item Soit $d\in\N$. Alors, d'après la propriété \ref{prop_lien_Euler_Bernoulli},
\[E_d(X+1)+E_d=\dfrac{2}{d+1}\left[B_{d+1}(X+1)+B_{d+1}(X)-2^{d+1}\left(B_{d+1}\left(\dfrac{X+1}{2}\right)+B_{d+1}\left(\dfrac{X}{2}\right)\right)\right]\]
Or, d'après la propriété \ref{prop_Bernoulli_un_demi}, $2^{d+1}\left(B_{d+1}\left(\frac{X+1}{2}\right)+B_{d+1}\left(\frac{X}{2}\right)\right)=2B_{d+1}$ donc $E_d(X+1)+E_d=\frac{2}{d+1}\left[B_{d+1}(X+1)-B_{d+1}\right]$ et, d'après la propriété \ref{prop_différence_Bernoulli},  $B_{d+1}(X+1)-B_{d+1}=(d+1)X^d$ donc $E_d(X+1)+E_d=2X^d$.
\item Soit $d\in\N^*$. Alors, en évaluant l'égalité précédente en $0$, il vient $E_d(1)+E_d(0)=0$ donc $E_d(1)=-E_d(0)$.
\item Soit $d\in\N^*$. En dérivant l'égalité de la propriété \ref{prop_lien_Euler_Bernoulli}, il vient
$E_d'=\frac{2}{d+1}\left[B_{d+1}'-2^dB_{d+1}'\left(\frac{X}{2}\right)\right]$. Or, d'après la propriété \ref{prop_dérivée_Bernoulli}, $B_{d+1}'=(d+1)B_{d}$ donc $E_d'=2\left[B_d-2^dB_d\left(\frac{X}{2}\right)\right]=dE_{d-1}$.
\item Soit $d\in\N$ et $a\in\R$. Par une récurrence immédiate, on déduit du point 3. que, pour tout $k\in\llbracket 0, d\rrbracket$, $E_d^{(k)}=\frac{d!}{(d-k)!}E_{d-k}$. Dès lors, d'après la formule de Taylor,
\[E_d(X)=\sum_{k=0}^d \frac{E_d^{(k)}(a)}{k!}(X-a)^k = \sum_{k=0}^d \frac{d!E_{d-k}(a)}{k!(d-k)!}(X-a)^k=\sum_{k=0}^d \dbinom{d}{k} E_{d-k}(a) (X-a)^k.\]
 \hfill $\square$
\end{enumerate}

\begin{corollaire} --- Pour tout $d\in\N$, $E_d(1-X)=(-1)^dE_d$. De plus, si $d$ est un entier pair strictement positif, $E_d(0)=E_d(1)=0$.
\label{coro_symétrie_Euler}
\end{corollaire}

\emph{Preuve}. --- On peut déduire ce corollaire, par un simple calcul, à partir des propriétés \ref{prop_Bernoulli_un_demi} et \ref{prop_lien_Euler_Bernoulli}. On peut également utiliser l'argument d'algèbre linéaire suivant. Considérons l'application $\varphi : P \mapsto P(X+1)+P$ définie sur $\R[X]$. Clairement, $\varphi$ est un endomorphisme du $\R-$espace vectoriel $\R[X]$ et, pour tout $n\in\N$, $\varphi(X^n)$ est un polynôme de degré $n$ donc la famille $(\varphi(X^n))_{n\in\N}$ est une famille de polynômes échelonnés en degré: c'est donc une base de $\R[X]$. Il s'ensuit que $\varphi$ est un isomorphisme de $\R[X]$. 

Soit $d\in\N$. Alors, d'après la propriété \ref{prop_polynômes_Euler},
\[\varphi(E_d(1-X))=E_d(-X)+E_d(1-X)=2(-X)^d=(-1)^d(2X^d)=(-1)^d\varphi(E_d)=\varphi((-1)^dE_d)\] 
donc, par injectivité de $\varphi$, $E_d(1-X)=(-1)^dE_d$. 

Supposons, à présent, $d$ pair et non nul. On déduit de ce qui précède que $E_d(1)=(-1)^dE_d(0)=E_d(0)$. Or, d'après le point 2. de la propriété \ref{prop_polynômes_Euler}, $E_d(1)=-E_d(0)$ donc on conclut que $E_d(0)=E_d(1)=0$.

 \hfill $\square$

\begin{definition} --- Soit $d\in\N$. Le nombre $e_d:=2^dE_d\left(\frac{1}{2}\right)$ est appelé le nombre d'Euler d'indice $d$.
\end{definition}

\begin{propriete} --- Pour tout $d\in\N$, $e_d\in\Z$ et, de plus, si $d$ est impair, alors $e_d=0$.
\end{propriete}

\emph{Preuve}. --- Soit $d\in\N$ un entier impair. D'après le corollaire \ref{coro_symétrie_Euler}, $E_d\left(\frac{1}{2}\right)=(-1)^dE_d\left(\frac{1}{2}\right)=-E_d\left(\frac{1}{2}\right)$ donc $E_d\left(\frac{1}{2}\right)=0$ et, par suite, $e_d=0$ et, en particulier, $e_d\in\Z$. 

Soit à présent un entier naturel $n$ non nul Alors, d'après le point 4. de la propriété \ref{prop_polynômes_Euler},
\begin{align*}
\sum_{k=0}^{2n} \dbinom{2n}{k}e_k&=\sum_{k=0}^{2n} \dbinom{2n}{k} 2^kE_k\left(\dfrac{1}{2}\right)=\sum_{k=0}^{2n} \dbinom{2n}{k} 2^{2n-k}E_{2n-k}\left(\dfrac{1}{2}\right) \\
&=2^{2n}\sum_{k=0}^{2n} \dbinom{2n}{k} E_{2n-k}\left(\dfrac{1}{2}\right)\left(1-\dfrac{1}{2}\right)^k=2^{2n}E_{2n}(1).
\end{align*}
Or, d'après le corollaire \ref{coro_symétrie_Euler}, $E_{2n}(1)=0$ donc $\sum\limits_{k=0}^{2n} \binom{2n}{k}e_k=0$. Comme, pour tout entier impair $k$, $e_k=0$, $\sum\limits_{k=0}^{n} \binom{2n}{2k}e_{2k}=0$. De plus, $e_0=E_0\left(\frac{1}{2}\right)=1$ donc, la suite $(e_{2n})$ est définie par la relation de récurrence, $e_0=1$ et, pour tout $n\in\N^*$, $e_{2n}=-\sum\limits_{k=0}^{n-1} \binom{2n}{2k}e_{2k}$ donc cette suite est à valeur entières. Ainsi, si $d$ est un entier pair, $e_d\in\Z$. \hfill $\square$

\begin{corollaire} --- Soit $d$ un entier pair. Alors,
\begin{enumerate}
\item $E_d\in\Z[X]$;
\item $e_d$ est un entier impair.
\end{enumerate}
\label{coro_prop_arithmétique_Euler}
\end{corollaire}

\emph{Preuve} \hspace{1cm}
\begin{enumerate}
\item Remarquons que, d'après le point 4. de la propriété \ref{prop_polynômes_Euler},
\[E_d=\sum_{k=0}^d \dbinom{d}{k} E_{d-k}\left(\dfrac{1}{2}\right)\left(X-\dfrac{1}{2}\right)^k=\dfrac{1}{2^d}\sum_{k=0}^d \dbinom{d}{k} 2^{d-k}E_{d-k}\left(\dfrac{1}{2}\right)(2X-1)^k\]
donc
\begin{equation}
E_d=\dfrac{1}{2^d}\sum_{k=0}^d \dbinom{d}{k} e_{d-k}(2X-1)^k
\label{eq_expression_Euler_1/2}
\end{equation}
Comme les nombres d'Euler sont entiers, on en déduit que $2^dE_d \in\Z[X]$. Ainsi, pour tout $k\in\llbracket 0, d \rrbracket$, le dénominateur $s_k$ du coefficient d'indice $k$ de $E_d$ est une puissance de $2$. Or, par définition, pour tout $k\in\llbracket 0, d\rrbracket$, le coefficient d'indice $k$ de $E_d$ est 
\[c_k:=\dbinom{d}{k}\dfrac{2(1-2^{d+1-k})}{d+1-k}b_{d+1-k}=\dfrac{1}{d+1}\dbinom{d+1}{k}(1-2^{d+1-k})(2b_{d+1-k}).\]
D'après le corollaire \ref{coro_vS_C}, la valuation $2-$adique du dénominateur de $b_{d+1-k}$ est au plus 1 donc le dénominateur de $2b_{d+1-k}$ est impair. De plus, $d+1$ est impair donc $s_k$ est impair. Comme $s_k$ est une puissance de $2$, on en conclut que $s_k=1$ et donc $c_k\in\Z$.
\item En conservant les notations de la question précédente, par définition,
\[e_d=2^d\sum_{k=0}^d c_k\left(\dfrac{1}{2}\right)^k=c_d+2\sum_{k=0}^{d-1} c_k2^{d-1-k}\]
Or, $c_d=1$ donc $e_d$ est impair. \hfill $\square$
\end{enumerate}

\subsection{Irrationalité des racines}

\begin{theoreme}[Brillhart, 1969] --- Soit $d\in\N^*$.
\begin{enumerate}
\item Si $d$ est pair alors les seules racines rationnelles de $E_d$ sont $0$ et $1$.
\item Si $d$ est impair alors l'unique racine rationnelle de $E_d$ est $\frac{1}{2}$.
\end{enumerate}
\end{theoreme}

\emph{Preuve} \hspace{1cm}
\begin{enumerate}
\item Supposons que $d$ est pair. Alors, d'après le corollaire \ref{coro_prop_arithmétique_Euler}, $E_d$ est un polynôme unitaire de $\Z[X]$. Ainsi, d'après le corollaire \ref{coro_entier_alg}, les racines rationnelles de $E_d$ sont des entiers. On sait, d'après le corollaire \ref{coro_symétrie_Euler}, que $0$ et $1$ sont racines de $E_d$. Soit un entier $m \geqslant 2$. Alors, d'après le point 1. de la propriété \ref{prop_polynômes_Euler}, $E_d(m)=2(m-1)^d-E_d(m-1)$ et, par une récurrence immédiate, puisque $E_d(1)=0$,
\[E_d(m)=2\left[(m-1)^d-(m-2)^d+(m-3)^d-(m-4)^d+\cdots+(-1)^m\right].\]
Or, pour tout entier $k\geqslant 0$, $(k+1)^d-k^d>0$ donc $E_d(m)>0$ et, en particulier, $E_d(m) \neq 0$. Enfin, si $m\leqslant -1$ alors, d'après le corollaire \ref{coro_symétrie_Euler}, $E_d(m)=(-1)^mE_d(1-m)$ et, comme $1-m \geqslant 2$, d'après le cas précédent, $E_d(1-m)\neq 0$ donc $E_d(m)\neq0$.
On conclut donc que les seules racines rationnelles de $E_d$ sont 0 et 1.
\item Supposons que $d$ est impair. D'après le corollaire \ref{coro_prop_arithmétique_Euler}, $e_d=0$ donc $E_d\left(\frac{1}{2}\right)=0$ et ainsi $\frac{1}{2}$ est bien racine de $E_d$. De plus, en utilisant l'égalité \eqref{eq_expression_Euler_1/2} et le fait que $e_{d}=0$, il vient
\[2^dE_d=\sum_{k=0}^{d} \dbinom{d}{k} e_{d-k}(2X-1)^k=(2X-1)P_d(2X-1)\]
en posant $P_d:=\sum\limits_{k=1}^{d} \binom{d}{k} e_{d-k}X^{k-1}$. Comme les nombres d'Euler sont entiers, $P_d\in \Z[X]$. De plus, $P_d$ est unitaire car $e_0=1$ et son coefficient constant est $de_{d-1}$ qui est impair d'après le corollaire \ref{coro_prop_arithmétique_Euler}. Ainsi, d'après la propriété \ref{prop_RRT}, si $P_d$ admet une racine rationnelle alors cette racine est un entier impair et, par suite, si $E_d$ admet une racine rationnelle autre que $\frac{1}{2}$ alors cette racine est un entier. Soit un entier $m\geqslant 1$. Comme dans le cas précédent, on a 
\[E_d(m)=2\left[(m-1)^d-(m-2)^d+(m-3)^d-(m-4)^d+\cdots+(-1)^m+(-1)^{m+1}E_d(1)\right].\]
Or,  d'après le point 2. de la propriété \ref{prop_polynômes_Euler}, $E_d(1)=-E_d(0)=\frac{2(2^{d+1}-1)}{d+1}b_{d+1}$. Comme $d+1$ est pair, le numérateur de $\frac{2(2^{d+1}-1)}{d+1}$ est impair. Or, d'après le corollaire \ref{coro_vS_C}, $2$ divise le dénominateur de $b_{d+1}$ donc le dénominateur de $E_d(1)$ est pair et, en particulier, $E_d(1)$ n'est pas entier. Par suite, on en déduit que $E_d(m)$ n'est pas entier et donc $E_d(m)\neq 0$. Par le même argument que dans le point 1., on en déduit que, si $m$ est un entier négatif ou nul, $E_d(m)=(-1)^mE_d(1-m)\neq 0$. 

On conclut donc que $\frac{1}{2}$ est l'unique racine rationnelle de $E_d$. \hfill $\square$
\end{enumerate}

\begin{rem} --- Comme pour les polynômes de Bernoulli, il est facile de voir que, pour tout entier $k\geqslant 1$, les seules racines de $E_{2k+1}$ dans $\intervalleff{0}{1}$ sont $0$ et $1$ et l'unique racine de $E_{2k}$ dans $\intervalleff{0}{1}$ est $\frac{1}{2}$. Le fait que ces nombres soient les seules racines rationnelles des polynômes d'Euler a été démontré par Brillhart dans son article \cite{Bri69}. Il a également prouvé qu'à l'exception de $E_5=\left(X-\frac{1}{2}\right)(X^2-X+1)^2$, les polynômes d'Euler n'ont que des racines simples. Une étude sur la localisation des racines réelles strictement supérieures à $1$ a été ménée par Delange \cite{Del88}.

Pour ce qui est de l'irréductibilité, on conjecture que, pour tout entier $k\geqslant 2$, le polynôme $\frac{E_{2k}}{X(X-1)}$ est irréductible sur $\Q$ et que, pour tout entier $k\geqslant 1$ différent de $2$, $\frac{E_{2k+1}}{X-\frac{1}{2}}$ est irréductible sur $\Q$. Des résultats partiels en ce sens ont été obtenus par Carlitz \cite{Car52} et Brillhart \cite{Bri69}, \cite{Bri72}.
\label{rem_Euler_irré}
\end{rem}

\bibliographystyle{smfalpha}
\bibliography{biblio} 

\end{document}